\documentclass[a4paper,UKenglish]{article}

\usepackage[utf8]{inputenc}
\usepackage[T1]{fontenc}
\usepackage{babel}
\usepackage{hyperref}

\usepackage{amsthm}
\usepackage{amsmath, amssymb}
\usepackage{cleveref}
\usepackage{bbm}
\usepackage{tikz}
\usepackage{algorithm}
\usepackage[noend]{algpseudocode}
\usepackage{thmtools}

\algrenewcommand\algorithmicrequire{\textbf{Input}}
\algrenewcommand\algorithmicensure{\textbf{Output}}

\definecolor{combi-cyan}{RGB}{0,170,170}
\definecolor{combi-darkcyan}{RGB}{0,100,100}
\definecolor{combi-orange}{RGB}{250,170,30}
\definecolor{combi-green}{RGB}{180,210,50}

\title{Parameterized complexity of scheduling unit-time jobs with generalized precedence constraints}
\author{
	Christina B\"using\thanks{RWTH Aachen University, Germany — \texttt{buesing@combi.rwth-aachen.de}} \and
	Maurice Draeger\thanks{RWTH Aachen University, Germany — \texttt{maurice.draeger@rwth-aachen.de}} \and
	Corinna Mathwieser\thanks{Corresponding author, RWTH Aachen University, Germany — \texttt{mathwieser@combi.rwth-aachen.de}}
}
\date{}

\theoremstyle{plain} 
\newtheorem{theorem}{Theorem}
\newtheorem{lemma}[theorem]{Lemma}

\newtheorem{corollary}[theorem]{Corollary}

\theoremstyle{definition} 
\newtheorem{definition}[theorem]{Definition}

\theoremstyle{remark} 

\begin{document}

\maketitle

\begin{abstract}
We study the parameterized complexity of scheduling unit-time jobs on parallel, identical machines under generalized precedence constraints for minimization of the makespan and the sum of completion times \linebreak($P \vert gen$-$prec, p_j=1\vert \gamma$, $\gamma\in\{C_{max},\sum_jC_j\}$). In our setting, each job is equipped with a Boolean formula (precedence constraint) over the set of jobs. A schedule satisfies a job's precedence constraint if setting earlier jobs to true satisfies the formula. Our definition generalizes several common types of precedence constraints: classical $and$-constraints if every formula is a conjunction, $or$-constraints if every formula is a disjunction, and $and/or$-constraints if every formula is in conjunctive normal form. We prove fixed-parameter tractability when parameterizing by the number of predecessors. For parameterization by the number of successors, however, the complexity depends on the structure of the precedence constraints. If every constraint is a conjunction or a disjunction, we prove the problem to be fixed-parameter tractable. For constraints in disjunctive normal form, we prove $W[1]$-hardness. We show that the $and/or$-constrained problem is $\mathcal{NP}$-hard, even for a single successor. Moreover, we prove $\mathcal{NP}$-hardness on two machines if every constraint is a conjunction or a disjunction. This result not only proves para-$\mathcal{NP}$-hardness for parameterization by the number of machines but also complements the polynomial-time solvability on two machines if every constraint is a conjunction~(\cite{coffman1972optimal}) or if every constraint is a disjunction (\cite{johannes2005complexity}).
\end{abstract}

\section{Introduction}
Scheduling jobs with precedence constraints on identical parallel machines is a well-researched problem with many real-world applications. A classical precedence constraint between two jobs, say from job $i$ to job $j$, means that job $j$ cannot begin until job $i$ has finished. If job $j$ has multiple such constraints, it can only start once all the predecessor jobs have been completed. 
However, many real-life settings find classical precedence constraints (\emph{and}-type constraints) too restrictive. It may not always be necessary to complete all predecessors before starting a job. This need for flexibility has led to the proposal of alternative variants of precedence constraints in the literature. The first variant is \emph{or}-precedence constraints where completing any predecessor is sufficient to start a job. The second variant includes \emph{and+or}-constraints where some jobs require all predecessors to be completed earlier ($and$-jobs) while others only require to be preceded by only one predecessor ($or$-jobs). Another more general variant is known as \emph{and/or}-constraints. In this model, each job has multiple groups of predecessor jobs, and to start the job, at least one job from each group must be completed first. A motivating example for $and/or$-constraints is provided by Goldwasser, Lotombe and Motwani (\cite{goldwasser1996complexity}), who examine the problem of disassembling a product. Certain parts may need to be taken out prior to accessing other components. Occasionally, the same part might be accessible from various directions and therefore, removing any of several different components may provide access to a particular part. Scheduling with a combination of $or$- and $and$-constraints has further been studied in the works of Gillies and Liu \cite{gillies1995scheduling}, M\"ohring et al. \cite{mohring2004scheduling}, Lee et al. \cite{lee2012flexible} and Erlebach et al. \cite{erlebach2003scheduling}. For scheduling under solely $or$-constraints, we refer to the works of Berit \cite{johannes2005complexity} and the paper by Happach \cite{happach2021makespan}.  \\
Given the diverse landscape of precedence constraint types, we propose a more general notion of precedence constraints, which encompasses all aforementioned types. Specifically, we encode precedence constraints regarding a specific job $ j $ as a logical formula $ \psi_j $ over (unnegated) variables representing the set $J$ of jobs: $ \{x_i \mid i \in J\} $. In any feasible schedule, $ \psi_j $ must be satisfied by the assignment resulting from setting $ x_i = 1 $ precisely if job $ i $ is executed prior to $ j $.
Our definition gives rise to $and$-constraints when each $ \psi_j $ is a conjunction; to $or$-constraints when each $ \psi_j$ is a disjunction; to \emph{and+or}-constraints when each $ \psi_j$ is either; and to $and/or$-constraints when each $ \psi_j$ takes conjunctive normal form. Naturally, another interesting variant arises when each $ \psi_j$ is expressed in disjunctive normal form. We refer to this new variant as $or/and$-constraints. To gain an intuitive motivation of $or/and$-constraints, consider the following example from disaster management: There are several routes that can be used to evacuate people from a village, but obstacles (such as debris or fallen trees) are blocking these routes. To proceed with the evacuation, it might not be essential to clear all routes yet but it is necessary to remove all obstacles from one route.\\
The study of parameterized algorithms as an approach to addressing computationally hard problems is a vibrant field of research that has gained significant momentum in the scheduling community in recent years (\cite{mnich2018parameterized}, \cite{prot2018survey}, \cite{chen2017parameterized}, \cite{bessy2019parameterized}, \cite{nederlof_et_al:LIPIcs.IPEC.2020.25}, \cite{mnich2015scheduling}, \cite{knop2018scheduling}).  In the context of classical precedence-constrained scheduling problems, the most common parameters include the number of machines, the width of the partial order that is induced by the precedence constraints, and the maximum processing time. Observe that the effect of precedence constraints on the complexity of scheduling problems can be twofold: On the one hand, the additional structure imposed by precedence constraints can simplify scheduling problems. An extreme example is when precedence constraints result in a total order among all jobs, rendering scheduling problems trivial. On the other hand, precedence constraints can also turn polynomial-time solvable scheduling problems into $\mathcal{NP}$-hard problems (Prot and Bellenguez-Morineau \cite{prot2018survey}): Scheduling unit-time jobs on parallel identical machines while minimizing the makespan or the (weighted) sum of completion times becomes $\mathcal{NP}$-hard in the presence of $and$-precedence constraints (Lenstra and Kan \cite{lenstra1978complexity}, Sethi \cite{sethi1977complexity}). 
Restricting the width of the partial order in classical precedence-constrained scheduling seeks to exploit the former effect: A small width generates high interrelatedness between the jobs and may thus provide scheduling problems with more structure. Yet, this approach does not necessarily result in fixed-parameter tractability, not even in combination with other classical parameters. Bodlaender and Fellows (\cite{bodlaender1995w}) show that even $P |prec, p_j =1|C_{max}$ is $W[2]$-hard when parameterized by the width and the number of machines and Van Bevern et al. show in \cite{van2016precedence} that $P2|prec, p_j \in \{1, 2\}|C_{max}$ is $W[2]$-hard when parameterized by the width. The diversity in the forms that precedence constraints can take—and the observation that enforcing dense dependency structures does not necessarily improve tractability—motivates the introduction of new parameters for generalized precedence constraints which limit interrelatedness.

	\subsection*{Our Contribution}
	We study the problem $P \vert gen$-$prec, p_j=1\vert \gamma$ of scheduling unit-time jobs on parallel, identical machines under generalized precedence constraints and consider the objectives of minimizing the makespan ($\gamma=C_{max}$) or the sum of completion times ($\gamma=\sum_jC_j$). We investigate the parameterized complexity with respect to three parameters: the number $k_p$ of predecessors, the number $k_s$ of successors and the number $m$ of machines. We achieve fixed-parameter tractability when parameterizing by the total number $k_p$ of jobs that appear as predecessors in any precedence constraint, by providing a single-exponential algorithm. For the number $k_s$ of successors however, the picture is more diverse and depends on the structure of the precedence constraints. For \emph{and+or}-constraints, we prove the problem to be fixed-parameter tractable (FPT) by providing an algorithm with doubly-exponential dependence on $k_s$. This result does not extend to $and/or$- nor to $or/and$-constraints: We prove that $P \vert and/or, p_j=1\vert \gamma$ is para-$\mathcal{NP}$-hard with respect to $k_s$ by showing that the problem is $\mathcal{NP}$-hard, even for a single successor. For $or/and$-constraints, we prove the problem to be in $XP$ and to be $W[1]$-hard. Since Berit provides in \cite{johannes2005complexity} a polynomial-time algorithm for scheduling unit-time jobs on parallel, identical machines under $or$-constraints, it makes sense to ask if \emph{and+or}-constrained scheduling remains FPT when only restricting the number $k_s^{and} \leq k_s$ of $and$-successors. However, we prove that scheduling under \emph{and+or}-constraints is $\mathcal{NP}$-hard, even if there is only one $and$-job. Since our reduction only makes use of two machines, we also obtain that $P \vert $\emph{and+or}$, p_j=1\vert \gamma$ is para-$\mathcal{NP}$-hard with respect to the number $m$ of machines. 

    \section{Preliminaries}
\subsection*{Problem Definition}
    We consider the problem of non-preemptively scheduling unit time jobs on identical, parallel machines under generalized precedence constraints. An instance specifies a set $ J = \{1, \ldots, n\} $ of jobs, a number $m$ of identical machines and occasionally job weights $ w_j \geq 0$, $j \in J$. A schedule $\mathcal{S}$ specifies which job, if any, each machine works on at each point in time. More precisely, a schedule $\mathcal{S}$ is represented as tuple $(C_j)_{j\in J}$, where $C_j \in \mathbb{N}$ indicates the \emph{completion time} of each job $j\in J$. Since $m$ jobs can be processed in parallel, we assume that $\vert \{j \in J\vert C_j=t\}\vert \leq m$ for any time $t\in \mathbb{N}$. We say that a job $i$ \emph{precedes} a job $j$ in $\mathcal{S}$ if $C_i < C_j$. Moreover, we refer to the period between time $t-1$ and time $t$ as \emph{time slot} $t\in \mathbb{N}$.\\
Precedence constraints comprehend the requirement that certain jobs have to be finished before others can be started. In order to formalize precedence constraints, we introduce a Boolean variable $x_j \in \{0,1\}$ for each job $j \in J$. Moreover, we denote by $\mathbbm{1}$ the logical one that represents the truth value ``true''. We encode precedence constraints as family $(\psi_j)_{j\in J}$, where each $\psi_j$, $j\in J$, is a Boolean formula over the (unnegated) variables in $\{x_i \vert i \in J\} \cup \{\mathbbm{1}\}$. If a variable $x_i$ occurs in a precedence constraint $\psi_j$, $i,j \in J$, then $i$ is called a \emph{predecessor of }$j$ and $j$ is called a \emph{successor of }$i$. A job $i \in J$, which is the predecessor of some other job, is called a \emph{predecessor}. Likewise, a job $j \in J$ with non-trivial precedence constraint $\psi_j\neq \mathbbm{1}$, is called a \emph{successor}. We refer by $k_p$ ($k_s$) to the number of predecessors (successors).\\
A schedule $\mathcal{S}$ is associated with variable assignments $x^j(\mathcal{S})$, $j\in J$, which set $x^j_i({\mathcal{S}})=1$ precisely if $i$ precedes $j$ in $\mathcal{S}$. A schedule $\mathcal{S}$ is called \emph{feasible} with respect to precedence constraints $(\psi_j)_{j\in J}$ if each $\psi_j$, $j \in J$, is satisfied by setting $x_i=x^j_i(\mathcal{S})$ for $i\in J$. Given a schedule $\mathcal{S}$, we call a job $j$ \emph{available} at time $t$ if the variable assignment that sets $x_i=1$ precisely if $C_i \leq t$ is a true assignment for $\psi_j$. We consider three classical objectives: Minimizing the makespan ($ \gamma:=C_{\max} = \max_{j \in J} C_j$), minimizing the sum of completion times ($ \gamma:=\sum_{j} C_j$) or minimizing the weighted sum of completion times, $ (\gamma:=\sum_{j} w_j C_j$). Following standard notation as introduced by Graham et al. in \cite{graham1979optimization}, we refer to the problem of finding an optimal schedule by $P \vert gen$-$prec ,p_j=1\vert \gamma$. Here, \( P \) means identical parallel machines, \( p_j = 1 \) indicates unit processing times for all jobs, and \( \gamma \in \{C_{\max},\, \sum_j C_j,\, \sum_j w_j C_j\} \) specifies the objective function. Given an instance $I$ of $P \vert gen$-$prec ,p_j=1\vert \gamma$, we denote the objective value of an optimal solution by $OPT(I)$ and of a specific schedule $\mathcal{S}$ by $\gamma(\mathcal{S})$. We distinguish five special cases of $P \vert gen$-$prec, p_j=1 \vert \gamma$, $\gamma \in \{C_{max}, \sum_jC_j, \sum_j w_j C_j\}$:

\begin{itemize}
	\item $P \vert prec  ,p_j=1\vert \gamma$, where each $\psi_j$ is a conjunction,
	\item $P \vert or$-$prec ,p_j=1 \vert \gamma$, where each $\psi_j$ is a disjunction,
	\item \emph{P}$\vert$\emph{and+or-}$prec  ,p_j=1\vert \gamma$, where each $\psi_j$ is either a conjunction or a disjunction,
	\item $P \vert and/or$-$prec ,p_j=1 \vert \gamma$, where each $\psi_j$ is in conjunctive normal form,
	\item $P \vert or/and$-$prec ,p_j=1 \vert \gamma$, where each $\psi_j$ is in disjunctive normal form.
\end{itemize}

\subsection*{Feasibility and perfect schedules}	
A feasible schedule for $P \vert gen$-$prec, p_j=1 \vert \gamma$ can be calculated in polynomial time by repeatedly choosing an available job and assigning it to any free machine. For a precise description, we refer to Algorithm 1 by M\"ohring et al. \cite{mohring2004scheduling}, which was formulated for $and/or$-constraints but can be easily implemented for generalized precedence constraints. In case of classical $and$-constraints, infeasibility is equivalent to the existence of a cycle in the precedence graph $G$, where $G$ is a digraph with node set $J$ and arc set $A=\{(i,j) \vert i,j \in J, \text{ $i$ is a predecessor of $j$}\}$. Throughout this paper, we assume that all instances are feasible. \\
We call a schedule $\mathcal{S}$ \emph{perfect} on $k \leq m$ machines if there exists a time step $t\in \mathbb{N}$, such that $k$ machines do not process a job after time $t$ but are never idle up to time $t-1$. We will make use of the following observation: If a schedule with optimum makespan $C_{max}$ is perfect on $m-1$ machines and processes $C_{max}$ many jobs on the remaining machine, then it simultaneously optimizes the (unweighted) sum of completion times. 

\subsection*{List Scheduling}	
List scheduling is a simple and intuitive method for addressing scheduling problems and was introduced by Graham in \cite{graham1966bounds}. The approach involves assigning a priority to each job, creating an ordered list based on these priorities. Jobs are then scheduled in the order they appear on this list. Whenever a machine is free, it is assigned the first unscheduled job from the list that is available to be processed at that time. List scheduling ensures machines are never left idle intentionally—idle time only occurs when no jobs are currently available for processing. This method can be efficiently implemented and runs in polynomial time.

\section{Parameterization by the number of \\predecessors}
		In this section, we analyze the parameterized complexity regarding the overall number $k_p$ of predecessors. 
        Note that the parameters $k_p$ and $k_s$ are, in a sense, conceptually opposed to width. While a small width typically indicates that most jobs are highly interrelated, having few predecessors or successors often leads to many jobs being largely independent. A small width is therefore used to enforce structure in scheduling problems that are hard without precedence constraints. In contrast, parameterizing by $k_p$ or $k_s$ aims to capture the effect of introducing a limited number of precedence constraints into problems that are efficiently solvable when no precedence constraints are present (as is the case for unit processing times).
        In the following, we prove that $P \vert$\emph{gen-prec}$, p_j=1 \vert \gamma$, $\gamma\in \{C_{max}, \sum_j C_j\}$ is FPT when parameterized by the number of predecessors (see Theorem \ref{thm: konstVorg}). For minimizing the weighted sum of completion times, we prove the weaker result that $P \vert$\emph{gen-prec}$, p_j=1 \vert  \sum_j w_jC_j$ lies in $\mathcal{P}$ if $k_p$ is bounded by a constant. In either case, we provide an algorithm that essentially explores all possibilities of positioning the predecessors in a schedule and later inserts all other jobs via List Scheduling. To test fewer possibilities, we argue that we may assume all predecessors to be processed before time step $k_p$ if $\gamma\in \{C_{max}, \sum_j C_j\}$. 
	
	\begin{lemma}\label{lem: predecessors_first}
		Let $I=(J, (\psi_j)_{j \in J})$ be a feasible instance of  $P \vert gen$-$prec, p_j=1 \vert \gamma$, $\gamma\in  \{C_{max}, \sum_j C_j\}$ with a set $J_p \subseteq J$ of $k_p$ predecessors. Then there exists an optimal schedule $\mathcal{S}^*$ with completion times $(C_j^*)_{j\in J}$ such that $C^*_j\leq k_p$ for all $j \in J_p$.
	\end{lemma}   
	
	\begin{proof}
		Let $\mathcal{S}^*$ be an optimal schedule with completion times $C_j^*$, $j \in J$. Assume $C_j^*>k_p$ for some $j \in J_p$. Let $t^*\leq k_p$ be the earliest time slot where no predecessor is processed. Choose a job $i \in J$ with $C_i^*=t^*$. By choice of $t^*$, $i$ is not a predecessor job. Let $h$ be an earliest predecessor with $C_h>t^*$. Construct a schedule $\mathcal{S}$ from $\mathcal{S}^*$ by swapping $i$ and $h$. Clearly, the swap is feasible and alters neither the makespan nor the sum of completion time. Thus $\mathcal{S}$ is optimal too, while the earliest time slot $t$ where no predecessor is processed has increased, $t>t^*$. Repeat until all completion times of predecessor jobs do not exceed $k_p$.
	\end{proof}

 	The algorithm {\sc alg-pre} for instances of $P \vert$\emph{gen-prec}$, p_j=1 \vert \gamma$, $\gamma \in  \{C_{max}, \linebreak \sum_j C_j, \allowbreak \sum_j w_j C_j\}$, with predecessors $J_p=\{i_1, \ldots, i_{k_p}\} \subseteq J$ proceeds as follows: Repeatedly pick a configuration $(C_{i_\ell})_{\ell=1}^{k_p}$ of predecessor completion times from $\mathcal{C}:=\{1,\ldots,k_p\}^{k_p}$ if $\gamma \in \{C_{max}, \sum_j C_j\}$ or from $\mathcal{C}:=\{1,\ldots,n\}^{k_p}$ if $\gamma=\sum_j w_jC_j$. If $(C_{i_\ell})_{\ell=1}^{k_p}$ is a feasible schedule for $J_p$, insert all remaining jobs through List Scheduling in order of non-increasing weight and return the schedule with the smallest objective value. A pseudocode description of the algorithm {\sc alg-pre} is depicted in Algorithm \ref{alg:konstVorg}. 
		
	\begin{restatable}{theorem}{konstVorg}\label{thm: konstVorg}
		Let $I$ be a feasible instance of $P \vert gen$-$prec, p_j=1 \vert \gamma$, $\gamma\in  \{C_{max},  \allowbreak \sum_j C_j, \sum_j w_j C_j\}$ with $k_p$ predecessors. If $\gamma\in  \{C_{max}, \sum_j C_j\}$, then {\sc alg-pre} returns in $\mathcal{O}(k_p^{k_p} \cdot m \cdot  n^2)$ time the completion times of an optimal schedule $\mathcal{S}^*$. Otherwise, {\sc alg-pre} returns in $\mathcal{O}(n^{k_p+2} \cdot m)$ time the completion times of an optimal schedule $\mathcal{S}^*$.
	\end{restatable}
	
	\begin{proof}
		Let $\mathcal{S}^*=(C_j^*)_{j\in J}$ be an optimal schedule, such that $(C^*_{i_1}, \ldots, \allowbreak C^*_{i_{k_p}}) \in \mathcal{C}$. (The existence of such $\mathcal{S}^*$ follows from Lemma \ref{lem: predecessors_first}). Consider the iteration of {\sc alg-pre} where $C_j=C_j^*$ for all $j \in J_p$. Clearly, the algorithm computes a feasible schedule $\mathcal{S}=(C_j)_{j \in J}$. Denote by $k_p^t$ the number of predecessor jobs processed during time slot $t\in \{1, \ldots, n\}$.\\
        Assume the jobs in $J \setminus J_p$ have $1\leq q\leq n-k$ different weights, that we denote by $w^1 \geq w^2 \geq \ldots \geq w^q$. We consider the vector $z(t)=(z_1(t), \ldots, z_q(t)) \in \mathbb{N}^q$ where $z_i(t)$, $i \in \{1, \ldots, q\}$, denotes the number of jobs in $J \setminus J_p$ with weight $w^i$ that $\mathcal{S}$ processes during time slot $t\in  \mathbb{N}$. For $\mathcal{S}^*$, we define $z^*(t)$, $t\in  \mathbb{N}$, analogously. If $z(t)=z^*(t)$ for all $t \in \mathbb{N}$, then $(C_{max}, \sum_j C_j, \sum_j w_jC_j)=(C^*_{max}, \sum_j C^*_j, \sum_j w_j C^*_j)$ and $\mathcal{S}$ is optimal. Otherwise, let $T:=\min\{t \vert z(t) \neq z^*(t)\}$ be the smallest time slot for which $z(T),z^*(T)$ are not equal. \\
		Assume first that $z(T)$ is lexicographically smaller than $z^*(T)$. Let $i\in \{1,\ldots,q\}$ be such that $z_s^*(T)=z_s(T)$ for $s<i$ and $z_i^*(T)>z_i(T)$. It holds that $\sum_{t=1}^T z_i(t)<\sum_{t=1}^Tz_i^*(t)$. Subsequently, there exists a job $j \in J$ with weight $w_j=w^i$, such that $C^*_j \leq T$ but $C_j > T$. Observe that $j$ is no predecessor, otherwise $C_j=C^*_j$. Since $C^*_j \leq T$, job $j$ is available at time $T-1$ in $\mathcal{S}^*$ and thus also in $\mathcal{S}$. Since {\sc alg-pre} sets $C_j > T$, the following holds: When $j$ is considered, the algorithm has already occupied all machines during time slot $T$, namely with predecessor jobs or jobs with weight at least $w_j=w^i$. This yields a contradiction as
		\begin{align*}
			m=k_p^T+\sum_sz_s(T)&=k_p^T+\sum_{s<i}z_s(T)+z_i(T)=k_p^T+\sum_{s<i}z^*_s(T)+z_i(T) \\&<k_p^T+\sum_{s<i}z^*_s(T)+z^*_i(T) \leq m. 
		\end{align*}
		We omit the analogous discussion of the case where $z^*(T)$ is lexicographically smaller than $z(T)$ (see \ref{app:konstVorg}) and conclude with a brief runtime analysis: We sort all non-predecessor jobs with respect to weight once ($\mathcal{O}(n\cdot log(n))$ time). Then the algorithm checks $\vert \mathcal{C} \vert$ configurations $(t_1,\ldots,t_{k_p})\in \mathcal{C}$ of predecessor completion times. For each such configuration, the algorithm inserts $\mathcal{O}(n)$ jobs $j\in J\setminus J_p$. For each job $j\in J\setminus J_p$, it checks for $\mathcal{O}(n)$ time steps $t$ whether $j$ is available at time $t$. For any such time step $t$, the algorithm searches for an idle machine among $m$ machines. Since, $\vert \mathcal{C} \vert \in \mathcal{O}(k^k)$ if $\gamma\in  \{C_{max}, \sum_j C_j\}$ and $\vert \mathcal{C} \vert \in \mathcal{O}(n^k)$ otherwise, the claimed asymptotic runtime follows.		
	\end{proof}
\section{Parameterization by the number of successors}
	We investigate whether restricting the number of successors has the same effect as restricting the number of predecessors. However, it turns out that the parameterized complexity with respect to the number of successors depends on the structure of the precedence constraints. More precisely, we prove that $P \vert $\emph{and+or}, $p_j=1 \vert\gamma$, $\gamma\in \{C_{max}, \sum_jC_j\}$ is FPT when parameterized with respect to the number $k_s$ of successors. For the more general $or/and$-constraints, we prove that $P \vert or/and, p_j=1 \vert \gamma$, $\gamma\in \{C_{max}, \sum_jC_j\}$ is polynomial-time solvable if the number $k_s$ of successors is bounded by a constant. We show that this result cannot be improved to fixed-parameter tractability by proving that $P \vert or/and, p_j=1 \vert \gamma$, $\gamma\in \{C_{max}, \sum_jC_j\}$ is $W[1]$-hard when parameterized by the number $k_s$ of successors. For $and/or$-constraints, we prove para-$\mathcal{NP}$-hardness. Throughout this section, we assume $\gamma \in \{C_{max}, \sum_j C_j\}$.\\
	
	To prove that \emph{P}$\vert$\emph{and+or-}$prec  ,p_j=1\vert \gamma$ is FPT with respect to $k_s$, we proceed as follows: 1. We provide an FPT algorithm {\sc and-alg-succ} for instances of $P \vert prec, p_j=1 \vert \gamma$ with $k_s$ successors. 2. We argue how to solve \emph{P}$\vert$\emph{and+or-}$prec  ,p_j=1\vert \gamma$ by repeatedly calling {\sc and-alg-succ}. Algorithm {\sc and-alg-succ} is similar to {\sc alg-pre} and explores all possibilities of positioning the successors in a schedule before inserting all other jobs. To avoid na\"{i}vely testing all $\mathcal{O}(n^{k_p})$ possible ways to schedule the successors, we use a similar result to Lemma \ref{lem: predecessors_first}: We may assume all successors to be processed within a time span of length $k_s$. 
	
	\begin{lemma}\label{lem: successors_span}
		Let $I=(J, (\psi_j)_{j \in J})$ be a feasible instance of $P \vert $\emph{gen-prec}$, p_j=1 \vert \gamma$, $\gamma\in  \{C_{max}, \sum_j C_j\}$ with $k_s$ successors and let $J_s\subseteq J$ denote the set of successors. Then there exists an optimal schedule $\mathcal{S}^*=(C^*_j)_{j \in J}$, such that $C^*_j-C^*_i\leq k_s-1$ for all $i,j \in J_s$.
	\end{lemma} 

	\begin{proof}
		Consider an optimal schedule $\mathcal{S}^*$ with completion times $\bar{C}_\ell$, $\ell\in J$, and assume that $\bar{C}^s_{max}-\bar{C}^s_{min}>k_s-1$ where $\bar{C}^s_{min}=\min\{C^*_\ell \vert \ell \in J_s\}$ is the completion time of an earliest successor and $\bar{C}^s_{max}=\max\{C^*_\ell \vert \ell \in J_s\}$ is the completion time of a latest successor. Let $t^*$ be the earliest time slot  with $\bar{C}^s_{min} < t^* < \bar{C}^s_{max}$ where no successor is processed. Pick a job $j \in J \setminus J_s$ with $C^*_j=t^*$. Let $i$ be a successor with $C^*_i=t^*-1$. We construct from $\mathcal{S}^*$ another feasible schedule $\mathcal{S}$ with completion times ${C}_\ell$, $\ell\in J$, by swapping jobs $i$ and $j$. Clearly, $\mathcal{S}$ is feasible and optimal and $C_{max}^s:=\max\{C_\ell \vert \ell \in J_s\}=\bar{C}^s_{max}$. If $t^*-1=\bar{C}^s_{min}$ and $i$ is the only successor with $\bar{C}_i=t^*_1-1$, then ${C}^s_{min}:=\min\{C_\ell \vert \ell \in J_s\}>\bar{C}^s_{min}$. Else, if $t^*-1>\bar{C}^s_{min}$ and $i$ is again the only successor with $\bar{C}_i=t^*_1-1$, then ${C}^s_{min}=\bar{C}^s_{min}$ but the earliest time slot $t\in ({C}^s_{min} , {C}^s_{max})$ without successors in $\mathcal{S}$ has decreased, $t<t^*$. Otherwise, we also get ${C}^s_{min}=\bar{C}^s_{min}$ while the overall number of time slots $t\in ({C}^s_{min} , {C}^s_{max})$ without successors in $\mathcal{S}$ has decreased. By repeating the argument, we obtain an optimal solution where all successors are processed within a time span of length $k_s$.
	\end{proof}

	Given a set $J_s \subseteq J$ of successors, the algorithm {\sc and-alg-succ} repeatedly picks a configuration $(C_j)_{j\in J_s}$ from 
    \begin{equation*}
     \mathcal{C}:=\{(t_j)_{j\in J_s} \in \{1, \ldots, n\}^{k_s}\vert \forall i,j \in J_s: t_i-t_j \leq k_s -1\},    
    \end{equation*}
    such that $\vert \{j \in J_s\vert C_j=t\}\vert  \leq m$ for all $t \in \{1, \ldots, n\}$ and $t_i<t_j$ if $i$ is a predecessor of $j$.  
    Then the algorithm inserts the remaining jobs as follows: {\sc and-alg-succ} iterates over all successors in order of non-decreasing completion times. For each successor $j \in J_s$, the algorithm assigns the earliest possible completion time to each predecessor. If some predecessor cannot be scheduled before $j$, the respective schedule of successors is discarded. In the end, all remaining jobs (which are neither successor nor predecessor) are inserted as early as possible. A pseudocode description of {\sc and-alg-succ} can be found in Algorithm \ref{alg:konstNach}.
		
	\begin{restatable}{proposition}{konstNachand}\label{thm: konstNach-and}
		Let $I$ be a feasible instance of $P \vert prec, p_j=1 \vert \gamma$, $\gamma\in  \{C_{max}, \allowbreak \sum_j C_j\}$ with $k_s$ successors. Then {\sc and-alg-succ} returns in $\mathcal{O}(k_s^{k_s} \cdot m \cdot  n^3)$ time an optimal schedule $\mathcal{S}^*$.
	\end{restatable}
\begin{proof}
Proposition \ref{thm: konstNach-and} can be shown along the same lines as Theorem \ref{thm: konstVorg}, see \ref{app:konstNachfolgand}.
\end{proof}

Using {\sc and-alg-succ}, it is easy to see that $P\vert or/and$-$prec, p_j=1 \vert \gamma$, $\gamma \in \{C_{max}, \sum_jC_j\}$, lies in $\mathcal{P}$ if the number of successors $k_s \leq const$ is bounded by some constant $const \geq 0$: Recall that $or/and$-constraints require each precedence constraint $\psi_j$, $j \in J$, to be a disjunction of conjunctive clauses. By electing one clause per job, we can solve an instance $I$ of $P\vert or/and$-$prec, p_j=1 \vert \gamma$ by solving  $\mathcal{O}(\vert I \vert ^{k_s})$ instances of $P\vert prec, p_j=1 \vert \gamma$ with at most $k_s$ successors. \\
For \emph{and+or}-constraints, we can strengthen this argument to prove that \linebreak \emph{P}$\vert$\emph{and+or-}$prec  ,p_j=1\vert \gamma$, $\gamma \in \{C_{max}, \sum_jC_j\}$ is FPT. Consider an instance $I=(J, (\psi_j)_{j\in J})$ of \emph{P}$\vert$\emph{and+or-}$prec  ,p_j=1\vert \gamma$ with $k_s$ successors. Let $J_s=J_s^{and} \dot \cup J_s^{or}$ be a partition of the set of successors, where a precedence constraint $\psi_j$ of a job $j \in J_s$ is a conjunction if $j \in J_s^{and}$ and a disjunction if $j \in J_s^{or}$. We call jobs in $J_s^{and}$ $and$-jobs and jobs in $J_s^{or}$ $or$-jobs. Consider an $or$-job $j \in J_s^{or}$ and let $i\in J$ be a predecessor of $j$. Consider the alternate instance $I'=(J, (\psi'_j)_{j\in J})$ which we obtain by fixing $i$ as only predecessor of $j$, that is, $\psi'_\ell=\psi_\ell$ if $\ell \neq j$ and $\psi'_j=x_i$. Clearly, any optimal schedule $\mathcal{S}'$ for $I'$ is feasible for $I$ too. However, whether $\mathcal{S}'$ is optimal for $I$ too, strongly depends on whether $i$ is also a predecessor of other jobs in $J_s\setminus \{j\}$. See Figure \ref{fig:constsuccand+or} for an example. Consider an arbitrary feasible schedule $\mathcal{S}$ for $I$. For each $or$-job $j \in J_s^{or}$, we can determine a predecessor job $i_j$ of $j$ such that $i_j$ precedes $j$ in $\mathcal{S}$. Observe that we can possibly have $i_j=i_k$ for two different $or$-jobs $j,k \in J_s^{or}$, $j\neq k$, and that each $i_j$ may be a predecessor of multiple $and$-successors. In other words, the choice of $(i_j)_{j \in J_s^{or}}$ gives rise to a function $\sigma:J_s^{or}\rightarrow2^{J_s}$ where $k \in \sigma(j)$ precisely if $k$ is an $or$-job and $i_j=i_k$ or if $k$ is an $and$-job with predecessor $i_j$. Our algorithm essentially tests all relevant variations of $\sigma$ and picks an $and$-predecessor $i_j \in J$ for each $or$-job $j \in J_s$ according to $\sigma$, before it applies {\sc and-alg-succ} to the resulting instance with $and$-constraints. 

\begin{figure}
    \centering
    \begin{tikzpicture}

        \filldraw[color=black, fill=blue!50, very thin] (0,0) -- (0,1) -- (1,1) -- (1,0);

        \filldraw[color=black, fill=blue!50, very thin] (0,1) -- (0,2) -- (1,2) -- (1,1);
        \filldraw[color=black, fill=blue!50, very thin] (1,1) -- (1,2) -- (2,2) -- (2,1);

        \filldraw[color=black, fill=blue!50, very thin] (0,2) -- (0,3) -- (1,3) -- (1,1);
        \filldraw[color=black, fill=orange!50, very thin] (1,2) -- (1,3) -- (2,3) -- (2,2);
        \filldraw[color=black, fill=teal!50, very thin] (2,2) -- (2,3) -- (3,3) -- (3,2);

        \path[-]{
            (0,0) edge (6,0)
            (0,1) edge (6,1)
            (0,2) edge (6,2)
            (0,3) edge (6,3)
            (0,0) edge (0,3)
        };   

        \node at (0.5, 2.5) {$a$};
        \node at (0.5, 1.5) {$b$};
        \node at (0.5, 0.5) {$c$};

        \node at (1.5, 2.5) {$d$};
        \node at (1.5, 1.5) {$e$};
        
        \node at (2.5, 2.5) {$f$};

        \filldraw[color=black, fill=blue!50, very thin] (7,0) -- (7,1) -- (8,1) -- (8,0);
        \filldraw[color=black, fill=blue!50, very thin] (7,1) -- (7,2) -- (8,2) -- (8,1);
        \filldraw[color=black, fill=orange!50, very thin] (7,2) -- (7,3) -- (8,3) -- (8,2);

        \filldraw[color=black, fill=teal!50, very thin] (8,0) -- (8,1) -- (9,1) -- (9,0);
        \filldraw[color=black, fill=blue!50, very thin] (8,1) -- (8,2) -- (9,2) -- (9,1);
        \filldraw[color=black, fill=blue!50, very thin] (8,2) -- (8,3) -- (9,3) -- (9,2);

        \path[-]{
            (7,0) edge (13,0)
            (7,1) edge (13,1)
            (7,2) edge (13,2)
            (7,3) edge (13,3)
            (7,0) edge (7,3)
        };

        \node at (7.5, 0.5) {$c$};
        \node at (7.5, 1.5) {$b$};
        \node at (7.5, 2.5) {$a$};

        \node at (8.5, 0.5) {$f$};
        \node at (8.5, 1.5) {$e$};
        \node at (8.5, 2.5) {$d$};

\end{tikzpicture} 
    \caption{We process six jobs $a,b,c,d,e,f$ with \emph{and+or}-precedence constraints $\psi_e= x_a\land x_b\land x_c$, $\psi_f=x_a \lor x_d$ on three machines. The left schedule is optimal for the $and$-variation with $\psi'_f=x_d$. The schedule on the right is optimal for the $and$-variation with $\psi'_f=x_a$ and for the original instance. 
    }
    \label{fig:constsuccand+or}
\end{figure}

\begin{definition}
	Let $I=(J, (\psi_j)_{j \in J})$ be an instance of $P\vert and$+$or$-$prec  ,p_j=1\vert \gamma$, $\gamma \in \{C_{max}, \sum_jC_j\}$ where $J_s=J_s^{or}\dot\cup J_s^{and}$ is the set of successors. For each $or$-job $j \in J_s^{or}$, let $i_j \in J$ be a predecessor of $j$.
	
	\begin{enumerate}
		\item Let $\sigma:J_s^{or}\rightarrow2^{J_s}$ be a function such that $k \in \sigma(j)$, $j \in J_s^{or}$, precisely if $k\in J_s^{or}$ and $i_j=i_k$ or if $k\in J_s^{and}$ with predecessor $i_j$. We call $\sigma$ the \emph{common predecessor function} of $I$ with respect to $(i_j)_{j \in J_s^{or}}$.
		\item If $\sigma:J_s^{or}\rightarrow2^{J_s}$ is a common predecessor function of $I$ with respect to $(i_j)_{j \in J_s^{or}}$, we call $(i_j)_{j \in J_s^{or}}$ \emph{generating representatives} of $\sigma$. 
		\item Let $I'=(J, (\psi'_j)_{j\in J})$ be the $P\vert prec, p_j=1 \vert \gamma$ instance with $\psi'_j=\psi_j$ if $j \notin J_s^{or}$ and $\psi_j=x_{i_j}$ if $j \in J_s^{or}$. We call $I'$ the \emph{$and$-variation} of $I$ with respect to $(i_j)_{j \in J_s^{or}}$. 
	\end{enumerate} 
\end{definition}

Observe the following properties of a common predecessor function $\sigma:J_s^{or}\rightarrow2^{J_s}$: First, $j \in \sigma(j)$ for all $j \in J_s^{or}$. If $j,k \in J_s^{or}$ and $k\in  \sigma(j)$, then $\sigma(j)=\sigma(k)$. If $\sigma$ is the common predecessor function of $I$ with respect to $(i_j)_{J_s^{or}}$, then $i_j$ is a predecessor of every $k \in \sigma(j)$, $j \in J_s^{or}$. Lastly, note that $\sigma$ may have different generating representatives $(i_j)_{J_s^{or}}$, $(i'_j)_{J_s^{or}}$, $(i_j)_{J_s^{or}}\neq(i'_j)_{J_s^{or}}$. The following lemma proves that for a given common predecessor function, the choice of generating representatives is irrelevant (under certain conditions). 
\begin{restatable}{lemma}{represent}\label{lem: representatives}
	Let $I=(J, (\psi_j)_{j \in J})$ be a $P\vert and$+$or$-$prec  ,p_j=1\vert \gamma$ instance, $\gamma \in \{C_{max}, \sum_jC_j\}$, with successor set $J_s=J_s^{or}\dot\cup J_s^{and}$. Let $\sigma, \sigma':J_s^{or}\rightarrow2^{J_s}$ be {common predecessor functions} of $I$ with respect to $(i_j)_{j \in J_s^{or}}, (i'_j)_{j \in J_s^{or}}$ respectively. Let $L, L'$ be the $and$-variations of $I$ with respect to $(i_j)_{j \in J_s^{or}}, (i'_j)_{j \in J_s^{or}}$ respectively. If $I_s:=\{j \in J_s^{or} \vert i_j \in J_s\}=\{j \in J_s^{or} \vert i'_j \in J_s\}$ and $i_j=i'_j$ for all $j \in I_s$, then $L$ is feasible precisely if $L'$ is feasible. In this case,  $OPT(L')= OPT(L)$ if moreover $\sigma=\sigma'$ and $i_j=i'_k$ whenever $i'_j=i_k$, $j,k \in J_s^{or}$.
\end{restatable}

\begin{proof}
	Assume that $L$ is infeasible and let $C=(j_1, \ldots, j_r)$ be a cycle in the precedence graph, $j_0:=j_r$. Consider a job $j_\ell$ in the cycle, $\ell \in \{1, \ldots, r\}$. If $j_\ell \in J_s^{and}$ is an $and$-job in $I$, then $j_{\ell-1}$ is also a predecessor of $j_\ell$ in $L'$. If $j_\ell \in J_s^{or}$, then $j_{\ell-1}=i_j$. Since $j_{\ell-1}$ is a successor in $L$, it is also a successor in $I$ and thus $j_{\ell-1}=i_j=i'_j$ is also a predecessor of $j_\ell$ in $L'$. Therefore, the precedence graph of $L'$ contains the same cycle $C$ and $L'$ is infeasible.\\
	Assume now that $L$ is feasible, $\sigma=\sigma'$ and $i_j=i'_k$ whenever $i'_j=i_k$, $j,k \in J_s^{or}$. Let $(C_j)_{j \in J}$ be the completion times of an optimal schedule $\mathcal{S}$. We define a schedule $\mathcal{S'}$ with completion times $(C'_j)_{j \in J}$ for $L'$. We denote the sets in $\{\sigma(j)\cap J_s^{or}\vert j \in J_s^{or}\}$ by $S_1, \ldots, S_r$. Observe that $(S_\ell)_{\ell=1}^r$ is a partition of $J_s^{or}$ and that $S_\ell =\sigma(j)$ precisely if $S_\ell=\sigma'(j)$, $j\in J_s^{or}, \ell =1, \ldots, r$. For each $\ell=1, \ldots, r$, let $i(\ell),i'(\ell) \in J_s^{or}$ be the common predecessor with $i(\ell)=i_k, i'(\ell)=i'_k$ for all $k\in S_\ell$. We obtain $\mathcal{S}'$ from $\mathcal{S}$ by swapping the positions of $i(\ell)$ and $i'(\ell)$ for all $\ell =1, \ldots,r$. 
	Observe that this is in fact the same as setting $C'_{i'_k}=C_{i_k}$ and $C'_{i_k}=C_{i'_k}$ for $k \in J_s^{or}$. All other jobs $j\in J$ remain as in $\mathcal{S}$, i.e. $C'_j=C_j$. Note that the schedule $\mathcal{S'}$ is well-defined, since we require $i'(\ell)=i(p)$ if $i(\ell)=i'(p)$ for $\ell, p =1, \ldots, r$. 
	Clearly, $\gamma(\mathcal{S}) =\gamma(\mathcal{S}') $. It remains to show that $\mathcal{S}'$ is feasible; for details we refer to \ref{app:konstNachfolgand+or}. 
\end{proof}

We now describe the algorithm {\sc and+or-alg-succ} for \emph{P}$\vert$\emph{and+or-}$prec  ,p_j=1\vert \gamma$, $\gamma \in \{C_{max}, \sum_jC_j\}$, in more detail. First, {\sc and+or-alg-succ} picks completion times $C_j$, $j \in J_s$, for the successors in the same way as {\sc and-alg-succ}. For every $or$-job $j \in J_s^{or}$ with a predecessor $i \in J_s$ with $C_i<C_j$, we set $i_j=i$. 
Afterwards, {\sc and+or-alg-succ} iterates over all $or$-jobs $j_1, \ldots, j_{q} \in J_s^{or}$ in order of non-decreasing completion time $C_{j_1}\leq \ldots \leq C_{j_{q}}$. If no $i_{j_i}$ has been chosen yet, we pick a subset $S_j \subseteq \{j_{i}, \ldots, j_q \}\cup J_s^{and}$ of successor jobs with $j\in S_j$. Then we pick a job $i \in J$ such that the successors of $i$ in $\{j_{i}, \ldots, j_q \}\cup J_s^{and}$ are precisely the vertices in $S_j$. We set $i_k=i$ for all $k\in S_j\cap J_s^{or}$. (If such a job $i$ does not exist, $S_j$ is an infeasible choice.) Finally, the algorithm {\sc and+or-alg-succ} proceeds along the same lines as {\sc and-alg-succ} for the $and$-variation of $I$ w.r.t. $(i_j)_{j\in J_s^{or}}$. 

	\begin{restatable}{theorem}{konstNachandplor}\label{thm: konstNach-and+or}
	Let $I$ be a feasible instance of $P\vert and$+$or$-$prec  ,p_j=1\vert \gamma$, $\gamma\in  \{C_{max}, \sum_j C_j\}$ with $k_s$ successors. Then {\sc and+or-alg-succ} returns in $\mathcal{O}(k_s^{k_s}  \cdot2^{k_s^2}\cdot m \cdot  n^3)$ time the completion times of an optimal schedule $\mathcal{S}^*$.
\end{restatable}

\begin{proof}
	 For an extended version of the proof, we refer to \ref{app:konstNachfolgand+or}. Let $J_s=J_s^{or} \dot \cup J_s^{and}$ be the set of successors with $\vert J_s \vert =k_s$ and $q:=\vert J_s^{or} \vert$. By Lemma \ref{lem: successors_span}, there exists an optimal schedule $\mathcal{S}^*=(C_j^*)_{j\in J}$ such that $(C^*_{i_1}, \ldots, C^*_{i_{k_p}}) \in \mathcal{C}$. Consider the iteration of {\sc and+or-alg-succ} where $C_j=C_j^*$ for all $j \in J_s$ and let $\mathcal{S}$ denote the schedule defined by $(C_j)_j$. Let $j_1, \ldots, j_q$ be the ordering, in which {\sc and+or-alg-succ} iterates over the $or$-jobs. We denote by $I_s$ the set of $or$-jobs $j \in J_s^{or}$, to which {\sc and+or-alg-succ} assigns a predecessor $i_j \in J_s$. \\
	The idea of the proof is the following: First, we pick $(i^*_j)_{j \in J_s^{or}}$ such that $i^*_j$ is a predecessor of $j$ with $C^*_{i^*_j} < C^*_j$, $j \in J_s^{or}$. We argue, that for a certain choice of $(S_j)_j$, the algorithm computes representatives $(i_j)_{j \in J_s^{or}}$ such that $(i_j)_j,(i^*_j)_j$ as well as the associated common predecessor functions $\sigma, \sigma^*$ comply with the conditions in Lemma \ref{lem: representatives}. \\
	More precisely, we pick $(i^*_j)_{j \in J_s^{or}}$ as follows: Consider the $or$-jobs $j$ in order $j=j_1, \ldots, j_q$. We set $i^*_j=i_j$ for $j \in I_s$. If $j=j_i$ and no $i^*_j$ has been chosen yet, we pick an arbitrary predecessor $i^*_j$ with $C^*_{i^*_j} < C^*_j$ but set $i^*_k=i^*_j$ for all $k=j_{i+1}, \ldots, j_q$, of which $i^*_j$ is also a predecessor. We denote by $\sigma^*:J_s^{or} \rightarrow 2^{J_s}$ the common predecessor function of $I$ w.r.t. $(i^*_j)_{j \in J_s^{or}}$ and consider the iteration where {\sc and+or-alg-succ} sets $S_{j_i}=\sigma^*(j_i)\setminus\{j_1, \ldots, j_{i-1}\}$ for $i=1, \ldots, q$. Observe that whenever {\sc and+or-alg-succ} is required to pick $i_j$ according to $S_j$, then $S_j=\sigma^*(j)$, otherwise $i_j$ would have been chosen in an earlier iteration. Due to the choice of $(S_j)_j$ and the properties of $(i^*_j)_j$, {\sc and+or-alg-succ} can find an $i_j$ for $j=j_r$ such that 1. $i_j$ is a predecessor of all jobs in $S_j$ and 2. $i_j$ is not a predecessor of any job in $(J_s^{and}\cup \{j_{r+1}, \ldots, j_q\})\setminus S_j$. We denote by $\sigma:J_s^{or} \rightarrow 2^{J_s}$ the common predecessor function of $I$ w.r.t. $(i_j)_{j \in J_s^{or}}$. Then it is easy to verify that $(i_j)_j,(i^*_j)_j,\sigma, \sigma^*$ comply with the conditions in Lemma \ref{lem: representatives}, which completes the proof.
	\end{proof}
We have shown that $P\vert or/and$-$prec, p_j=1 \vert \gamma$ is polynomial-time solvable if the number $k_s$ of successors is bounded by a constant while \emph{P}$\vert$\emph{and+or-}$prec,p_j=1\vert \gamma$ is even FPT with respect to the number of successors. Recall that Berit \cite{johannes2005complexity} proved the $or$-constrained problem $P\vert or, p_j=1 \vert \gamma$ to be polynomial-time solvable, even if no restrictions are made on the number of successors. Naturally, three questions arise:

\begin{enumerate}
	\item Is $P\vert or/and$-$prec, p_j=1 \vert \gamma$ FPT when parameterized by $k_s$? We give a negative answer to this question by proving that $P\vert or/and$-$prec, p_j=1 \vert \gamma$ is $W[1]$-hard (see Theorem \ref{thm:w1}). 
	\item Is $P\vert and/or$-$prec, p_j=1 \vert \gamma$ polynomial-time solvable if the number $k_s$ of successors is bounded by a constant? We will give a negative answer to this question by proving that $P\vert and/or$-$prec, p_j=1 \vert \gamma$ is $\mathcal{NP}$-hard, even if $k_s=1$ (see Proposition \ref{prop:complexity_andor_onesuccessor}).
    \item Is \emph{P}$\vert$\emph{and+or-}$prec  ,p_j=1\vert \gamma$ FPT with respect to the number of $and$-successors? Surprisingly, $P2\vert$\emph{and+or-}$prec  ,p_j=1\vert \gamma$ turns out to be para-$\mathcal{NP}$-hard when parameterized by the number of $and$-successors. In fact, we prove that the problem is $\mathcal{NP}$-hard, even if there is only a single $and$-successor (see Theorem \ref{thm:complexity_two_machines}). 
	
\end{enumerate}

\begin{restatable}{theorem}{weins}\label{thm:w1}
	$P\vert or/and$-$prec, p_j=1 \vert \gamma$, $\gamma \in \{C_{max}, \sum_jC_j\}$ is $W[1]$-hard when parameterized by the number $k_s$ of successors.
\end{restatable}
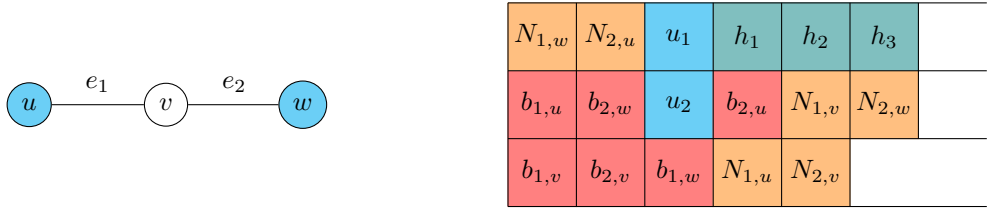
\begin{figure}[] \begin{center}  \begin{tikzpicture}[scale=.9] 

        \node[circle, draw, fill=cyan!50] (1) at (0,0.5) {$u$};
        \node[circle, draw] (2) at (2,0.5) {$v$};
        \node[circle, draw, fill=cyan!50] (3) at (4,0.5) {$w$};

        \path[-]{
        (1) edge (2)
        (3) edge (2)
        };

        \node at (1,0.8) {$e_1$};
        \node at (3,0.8) {$e_2$};

    \filldraw[red!50] (7,0) -- (7,-1) -- (8,-1) -- (8,0);
    \filldraw[red!50] (7,1) -- (7,0) -- (8,0) -- (8,1);
    \filldraw[orange!50] (7,1) -- (7,2) -- (8,2) -- (8,1);

    \filldraw[red!50] (8,0) -- (8,-1) -- (9,-1) -- (9,0);
    \filldraw[red!50] (8,1) -- (8,0) -- (9,0) -- (9,1);
    \filldraw[orange!50] (8,1) -- (8,2) -- (9,2) -- (9,1);

    \filldraw[red!50] (9,0) -- (9,-1) -- (10,-1) -- (10,0);
    \filldraw[cyan!50] (9,1) -- (9,0) -- (10,0) -- (10,1);
    \filldraw[cyan!50] (9,1) -- (9,2) -- (10,2) -- (10,1);
    
    \filldraw[orange!50] (10,0) -- (10,-1) -- (11,-1) -- (11,0);
    \filldraw[red!50] (10,1) -- (10,0) -- (11,0) -- (11,1);
    \filldraw[teal!50] (10,1) -- (10,2) -- (11,2) -- (11,1);
    
    \filldraw[orange!50] (11,0) -- (11,-1) -- (12,-1) -- (12,0);
    \filldraw[orange!50] (11,1) -- (11,0) -- (12,0) -- (12,1);
    \filldraw[teal!50] (11,1) -- (11,2) -- (12,2) -- (12,1);
    
    \filldraw[orange!50] (12,1) -- (12,0) -- (13,0) -- (13,1);
    \filldraw[teal!50] (12,1) -- (12,2) -- (13,2) -- (13,1);
    
    \path[-]{
    (7,-1) edge (14,-1)
    (7,0) edge (14,0)
    (7,1) edge (14,1)
    (7,2) edge (14,2)
    (7,-1) edge (7,2)
    (8,-1) edge (8,2)
    (9,-1) edge (9,2)
    (10,-1) edge (10,2)
    (11,-1) edge (11,2)
    (12,-1) edge (12,2)
    (13,0) edge (13,2)
    
    
    };

    \node at (7.5, 1.5){$N_{1,w}$};
    \node at (7.5, .5){$b_{1,u}$};
    \node at (7.5, -.5){$b_{1,v}$};
    
    \node at (8.5, 1.5){$N_{2,u}$};
    \node at (8.5, .5){$b_{2,w}$};
    \node at (8.5, -.5){$b_{2,v}$};

    \node at (9.5, 1.5){$u_1$};
    \node at (9.5, .5){$u_2$};
    \node at (9.5, -.5){$b_{1,w}$};

    \node at (10.5, 1.5){$h_1$};
    \node at (10.5, .5){$b_{2,u}$};
    \node at (10.5, -.5){$N_{1,u}$};

    \node at (11.5, 1.5){$h_2$};
    \node at (11.5, .5){$N_{1,v}$};
    \node at (11.5, -.5){$N_{2,v}$};

    \node at (12.5, 1.5){$h_3$};
    \node at (12.5, .5){$N_{2,w}$};

\end{tikzpicture}\end{center}
\caption{Reduction from Independent Set with $k=2$ to $P\vert or/and$-$prec, p_j=1 \vert \gamma$ with $k_s=2k+1$. Blue jobs and preceding orange jobs represent the size and selection of an independent set. In schedules with $C_{max}=2k+1$, green (red) jobs are used to block space after (before) blue jobs.}\label{fig:w_1}
\end{figure}

\begin{proof}
	An extended version of the proof is provided in \ref{app:konstNachfolgerorand}. We prove the claim by an FPT-reduction from Independent Set. Let $I'=(G,k)$ be an instance of Independent Set where $G=(V,E)$ is an undirected graph with $n':=\vert V \vert$ nodes and $k \in \mathbb{N}_0$ is the required size of an independent set, w.l.o.g. $k<n$. By $\overline{\mathcal{N}}(v):=\{w \in V \vert \{v,w\} \in E\} \cup \{v\}$ we denote the closed neighborhood of a vertex $v \in V$. We construct from $I'$ an instance $I$ of $P\vert or/and$-$prec, p_j=1 \vert C_{max}$ with $n'$ machines and require a makespan of at most $2(k+1)$. The construction for $\gamma=\sum_j C_j$ is analogous. An interpretation of the construction is given below.
    We introduce $2(k+k\cdot n')+1$ jobs by setting 
    \begin{align*}
        J:=\{u_i \vert i=1, \ldots, k\} \cup \{N_{iv},b_{iv} \vert i=1, \ldots, k, v\in V\} \cup \{h_i\vert i=1, \ldots, k+1\} 
    \end{align*}
	The precedence constraints are defined as follows: The auxiliary jobs $h_{1}, \ldots, h_{k+1}$ have to be scheduled in a chain that follows the completion of all $u_i$, $i=1, \ldots, k$. That is, $\psi_{h_{i+1}}=x_{h_i}$ for all $i=1, \ldots,k$ and $\psi_{h_1}=\bigwedge_{i=1}^k x_{u_i}$. The precedence constraint of a job $u_i$ is the disjunction $\psi_{u_i}=\bigvee_{v\in V}K_{iv}$ of $n'$ clauses $K_{iv}$, where 
    \begin{equation*}
        K_{iv}:=\bigwedge_{u \notin \overline{\mathcal{N}}(v)} x_{N_{iu}} \wedge \bigwedge_{u \in \overline{\mathcal{N}}(v)} x_{b_{iu}} \wedge \bigwedge_{j< i}x_{N_{jv}},\text{ }i=1, \ldots, k, v\in V.
    \end{equation*}
    We interpret $v\in V$ as the $i$-th node in the independent set if $K_{iv}$ is satisfied by the jobs that precede $u_i$, $i=1, \ldots, k$. In a schedule with makespan $2(k+1)$, the chain of the auxiliary jobs $h_1,  \ldots, h_{k+1}$ guarantees that all jobs $u_1, \ldots, u_k$ have to be completed by time $k+1$. This leaves a total of at most $k \cdot n'$ time slots to schedule the predecessors that are needed to make $u_1, \ldots, u_k$ available. However, each clause $K_{iv}$ is composed of $n'+i-1$ literals in a way that each of the last $i-1$ predecessors $N_{jv}$, $j=1, \ldots,i-1$, has to be a common predecessor with $u_j$, due to the space restrictions. Due to this, a node $v_i$ for which $K_{iv_i}$ is satisfied is disjoint and non-adjacent to the nodes $v_j$ for which the previous $i-1$ clauses $K_{jv_j}$ are satisfied, $j=1, \ldots, i-1$. Figure \ref{fig:w_1} depicts an example of the construction.
\end{proof}
From the proof of Theorem \ref{thm:w1}, we moreover get the following corollary: 

\begin{corollary}
    $P\vert or/and$-$prec, p_j=1 \vert C_{max}$ is $W[1]$-hard if parameterized by $C_{max}$ and $k_s$.
\end{corollary}

Given the complexity result in Theorem \ref{thm:w1}, the existence of a fixed-parameter tractable algorithm for $P\vert or/and$-$prec, p_j=1 \vert \gamma$ parameterized by $k_s$ is unlikely. In the context of $and/or$-constraints, we provide even stronger evidence against the fixed-parameter tractability of $P \vert and/or, p_j=1 \vert \gamma$ parameterized by $k_s$, by proving that the problem is para-$\mathcal{NP}$-hard.

\begin{restatable}{proposition}{konstNachandor}\label{prop:complexity_andor_onesuccessor}
    $P \vert and/or, p_j=1 \vert \gamma$, $\gamma \in \{C_{max}, \sum_j C_j\}$ is $\mathcal{NP}$-hard, even if $k_s=1$.
\end{restatable}
		
\begin{proof}
    We prove the claim by a reduction from Vertex Cover.
    Let $I'=(G,k)$ be an instance of Vertex Cover, where $G=(V,E)$ is an undirected Graph and $k \in \mathbb{N}$ is the required size of a vertex cover. We may assume that $k<|V|=:n'$, otherwise $I$ is trivial.
    We construct from $I'$ an instance $I$ of $P \vert and/or, p_j=1 \vert  \gamma$ with a single successor.
    We consider $m:=n'$ machines. The set $J$ of jobs is composed of a job $v$ for each vertex $v\in V$, a single job $e$, which will be the only successor and $n'-k$ auxiliary jobs ${b_1,..., b_{n'-k}}$. Through precedence constraints, we force all auxiliary jobs as well as one end-vertex per edge to be processed before $e$: We set $\psi_e=\bigwedge_{i=1}^{n'-k}x_{b_i} \land \bigwedge_{\{v,w\}\in E} (x_v \lor x_w)$ and $\psi_j=\mathbbm{1}$ for all $j\in J\setminus \{e\}$.
    We require to find a feasible schedule $\mathcal{S}$ with $C_{max}\leq 2$ if $\gamma=C_{max}$ or $\sum_{j\in J} C_j\leq 3n'-2k+2$ if $\gamma=\sum_j C_j$.
    It is easy to see that $I'$ is a yes-instance precisely if we can find a schedule for $I'$, which processes all auxiliary jobs and a vertex cover of size $k$ during the first time slot and all remaining jobs during time slot two. The construction is visualized in Figure \ref{fig:constsuccand-or}.    
\end{proof}
From the proof of Proposition \ref{prop:complexity_andor_onesuccessor}, we also get the following corollary:
\begin{figure}
    \centering
    \begin{tikzpicture}

        \node[circle, draw, fill=cyan!50] (1) at (0,0.5) {$v_1$};
        \node[circle, draw, fill=orange!50] (2) at (2,0.5) {$v_2$};
        \node[circle, draw, fill=cyan!50] (3) at (4,0.5) {$v_3$};

        \path[-]{
        (1) edge (2)
        (3) edge (2)
        };

        \node at (1,0.8) {$e_1$};
        \node at (3,0.8) {$e_2$};

    \filldraw[red!50] (7,0) -- (7,-1) -- (8,-1) -- (8,0);
    \filldraw[red!50] (7,1) -- (7,0) -- (8,0) -- (8,1);
    \filldraw[orange!50] (7,1) -- (7,2) -- (8,2) -- (8,1);

    \filldraw[cyan!50] (8,0) -- (8,-1) -- (9,-1) -- (9,0);
    \filldraw[cyan!50] (8,1) -- (8,0) -- (9,0) -- (9,1);
    \filldraw[teal!50] (8,1) -- (8,2) -- (9,2) -- (9,1);
    
    \path[-]{
    (7,-1) edge (12,-1)
    (7,0) edge (12,0)
    (7,1) edge (12,1)
    (7,2) edge (12,2)
    (7,-1) edge (7,2)
    
    };

    \node at (7.5, 1.5){$v_2$};
    \node at (7.5, .5){$b_1$};
    \node at (7.5, -.5){$b_2$};
    \node at (8.5, 1.5){$e$};
    \node at (8.5, .5){$v_1$};
    \node at (8.5, -.5){$v_3$};

\end{tikzpicture}

    \caption{Reduction of Vertex Cover with $n'=3$ and $k=1$ to $P|and/or$-$prec, p_j=1|\gamma$ with $k_s=1$}
    \label{fig:constsuccand-or}
\end{figure}
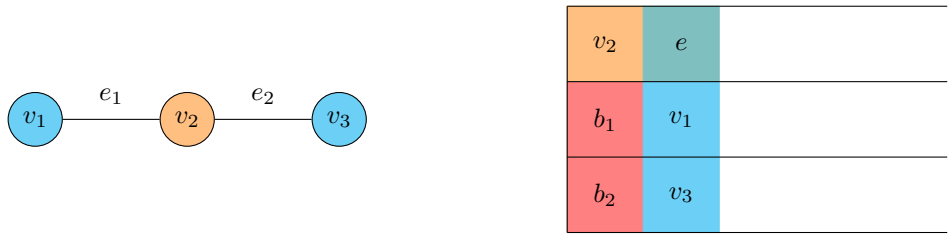
\begin{corollary}
    $P \vert and/or, p_j=1 \vert C_{max}$ is para-$\mathcal{NP}$-hard, if parameterized by $C_{max}$ and $k_s$.
\end{corollary}

\section{Parameterization by the number of machines}

	A classical parameter in the study of parameterized complexity of scheduling problems, is the number of machines. In \cite{coffman1972optimal}, Coffman and Graham prove that the $\mathcal{NP}$-hard problem $P \vert prec, p_j=1 \vert C_{max}$ becomes polynomial when restricted to two machines. For an arbitrary but constant number of machines, the complexity is open. However, Bodlaender and Fellows \cite{bodlaender1995w} prove that $P \vert and, p_j=1 \vert C_{max}$ is $W[2]$-hard when parameterized by the number of machines. In the following, we prove that $P2 \vert $\emph{and+or-prec}$, p_j=1 \vert \gamma$ is  $\mathcal{NP}$-hard for $\gamma=C_{max}, \sum_jC_j, \sum_jw_jC_j$ and thus $P2 \vert $\emph{and+or-prec}$, p_j=1 \vert \gamma$ is para-$\mathcal{NP}$-hard when parameterized by the number of machines. This result is also interesting in the context of the interplay between $and$- and $or$-constraints: On the one hand, it shows that the introduction of $or$-constraints turns the polynomial-time solvable $P2 \vert prec, p_j=1 \vert C_{max}$ into an $\mathcal{NP}$-hard problem. On the other hand, the result below shows that the introduction of a single $and$-job turns the polynomial-time solvable $P \vert or$-$prec, p_j=1 \vert C_{max}$ into an $\mathcal{NP}$-hard problem. For an \emph{and-or}-constrained instance with successors in $J_s=J_s^{or} \dot \cup J_s^{and}$, we denote by $k_s^{and}:=\vert J_s^{and} \vert$ the number of $and$-successors.
	
	\begin{restatable}{theorem}{twomachines}\label{thm:complexity_two_machines}
		$P2 \vert and$+$or$-$prec, p_j=1 \vert \gamma$, $\gamma\in \{C_{max}, \sum_j C_j\}$ with $k_s^{and}=1$ is $\mathcal{NP}$-hard.
	\end{restatable}    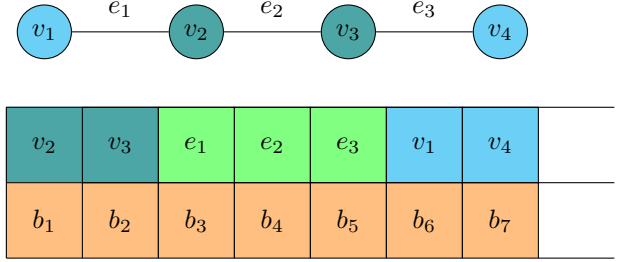
\begin{figure}[ht]
    \begin{center}
    \begin{tikzpicture}

        \node[circle, draw, fill=cyan!50] (1) at (.5,3) {$v_1$};
        \node[circle, draw, fill=teal!70] (2) at (2.5,3) {$v_2$};
        \node[circle, draw, fill=teal!70] (3) at (4.5,3) {$v_3$};
        \node[circle, draw, fill=cyan!50] (4) at (6.5,3) {$v_4$};

        \path[-]{
        (1) edge (2)
        (3) edge (2)
        (3) edge (4)
      
        };

        \node at (1.5,3.3) {$e_1$};
        \node at (3.5,3.3) {$e_2$};
        \node at (5.5,3.3) {$e_3$};

        \filldraw[teal!70] (0,1) -- (0,2) -- (1,2) -- (1,1);
        \filldraw[teal!70, draw=black] (1,1) -- (1,2) -- (2,2) -- (2,1);
        \filldraw[green!50, draw=black] (2,1) -- (2,2) -- (3,2) -- (3,1);
        \filldraw[green!50, draw=black] (3,1) -- (3,2) -- (4,2) -- (4,1);
        \filldraw[green!50, draw=black] (4,1) -- (4,2) -- (5,2) -- (5,1);
        \filldraw[cyan!50, draw=black] (5,1) -- (5,2) -- (6,2) -- (6,1);
        \filldraw[cyan!50, draw=black] (6,1) -- (6,2) -- (7,2) -- (7,1);

               \filldraw[orange!50] (0,0) -- (0,1) -- (1,1) -- (1,0);
        \filldraw[orange!50, draw=black] (1,0) -- (1,1) -- (2,1) -- (2,0);
        \filldraw[orange!50, draw=black] (2,0) -- (2,1) -- (3,1) -- (3,0);
        \filldraw[orange!50, draw=black] (3,0) -- (3,1) -- (4,1) -- (4,0);
        \filldraw[orange!50, draw=black] (4,0) -- (4,1) -- (5,1) -- (5,0);
        \filldraw[orange!50, draw=black] (5,0) -- (5,1) -- (6,1) -- (6,0);
        \filldraw[orange!50, draw=black] (6,0) -- (6,1) -- (7,1) -- (7,0);
    
        \path[-]{
            (0,0) edge (8,0)
            (0,1) edge (8,1)
            (0,2) edge (8,2)

            (0,0) edge (0,2)
        };

        \node at (0.5, 1.5) {$v_2$};
        \node at (1.5, 1.5) {$v_3$};
        \node at (2.5, 1.5) {$e_1$};
        \node at (3.5, 1.5) {$e_2$};
        \node at (4.5, 1.5) {$e_3$};
        \node at (5.5, 1.5) {$v_1$};
        \node at (6.5, 1.5) {$v_4$};

        \node at (0.5, 0.5) {$b_1$};
        \node at (1.5, 0.5) {$b_2$};
        \node at (2.5, 0.5) {$b_3$};
        \node at (3.5, 0.5) {$b_4$};
        \node at (4.5, 0.5) {$b_5$};
        \node at (5.5, 0.5) {$b_6$};
        \node at (6.5, 0.5) {$b_7$};
        
\end{tikzpicture}
    
\end{center}
    \caption{Reduction from Vertex Cover to $P2|$\emph{and+or-prec}$, p_j=1|\gamma$ with $J_s^{and}=\{b_6\}$.}
    \label{fig:constMachines}
\end{figure}
\begin{proof} 
	We prove the claim by a reduction from Vertex Cover. Let $I=(G,k)$ be an instance of Vertex Cover, where $G=(V,E)$ is an undirected graph and $k$ is the requested size of a vertex cover. We assume that $k<\vert V \vert$, otherwise $I$ is trivial. We construct from $I$ an instance $I'$ of $P2 \vert $\emph{and+or-prec}$, p_j=1 \vert \gamma$, $\gamma\in \{C_{max}, \sum_j C_j\}$, where the idea is the following: We introduce one job per edge, one job per vertex and $\ell:=\vert V \vert +  \vert E \vert$ auxiliary jobs to block the second machine. Using precedence constraints, we will force a feasible schedule to process all jobs in a vertex cover of size at most $k$ before finishing the edge jobs. 
    \\More precisely, we set $J:= V \cup E \cup B$ where $B=\{b_1, \ldots, b_{\ell}\}$ is the set of auxiliary jobs. All vertex jobs, as well as the first auxiliary job, may start at any time, i.e. $\psi_j=\mathbbm{1}$ if $j \in V \cup \{b_1\}$. Edges may not be scheduled unless one of their end-vertices has been scheduled already, i.e. $\psi_{e}=x_{u} \lor x_v$ for $e =\{u,v\} \in E$. We set $\psi_{b_{k+|E|+1}}=x_{b_{k+|E|}} \land \bigwedge_{e\in E} x_{e}$ to make sure that all edge jobs have been scheduled by the time we start the $(k+\vert E \vert+1)$-st auxiliary job. To process all auxiliary jobs in a chain, we set $\psi_{b_{i}}=x_{b_{i-1}} $, $i \in \{2, \ldots, \ell\} \setminus \{k+\vert E\vert+1\}$. Observe that $J_s^{and}:=\{b_{k+|E|+1}\}$ and $J_s^{or}:=E \cup \{b_1, \ldots b_{k+|E|}, b_{k+|E|+2}, \ldots, b_\ell\}$ is a partition of the successors into $or$- and $and$-jobs. The construction is visualized in Figure \ref{fig:constMachines}.
	Observe that a schedule is perfect (on all machines) precisely if it yields $C_{max}=\ell$ or, equivalently, if $\sum_j C_j= \ell\cdot(\ell+1)$. It is easy to see that $I$ is a yes-instance, precisely if there exists a schedule for $I'$ with $C_{max}=\ell$, ($\sum_j C_j= \ell\cdot(\ell+1)$); for details see \ref{app:twomachines}. 
\end{proof}

 \section{Summary and Outlook}
 In this paper, we study the problems of minimizing the makespan and the sum of completion times in parallel machine scheduling of unit-time jobs under precedence constraints. In our setting, each job $j$ is associated with a precedence constraint $\psi_j$, encoded as a Boolean formula over the set of jobs. A schedule satisfies $\psi_j$ if setting exactly the jobs scheduled before $j$ to true satisfies the formula. Depending on the structure of these formulas, we distinguish between $and$-, $or$-, \emph{and+or}-, $or/and$-, and $and/or$-constraints.\\
To analyze the parameterized complexity of these problems, we introduce two new structural parameters: the total number of predecessors and the total number of successors. Restricting the number of predecessors yields fixed-parameter tractability, even for general precedence constraints. In contrast, restricting the number of successors leads to a more nuanced landscape: for traditional $and$-constraints and more general \emph{and+or}-constraints, the problem is FPT, but for $or/and$-constraints (DNF), it becomes $W[1]$-hard. For and/or-constraints (CNF), the problem is even para-$\mathcal{NP}$-hard. A common parameter in the analysis of scheduling problems is the number of machines. While minimizing the makespan for $and$-constrained unit jobs is known to be $W[2]$-hard parameterized by the number of machines (Bodlaender and Fellows \cite{bodlaender1995w}), we show that the same problem is para-$\mathcal{NP}$-hard in the presence of \emph{and+or}-constraints.\\
We believe these results open up several promising research directions. One such direction concerns the interplay between $and$- and $or$-constraints. Naturally, any known hardness result for problems with classical $and$-constraints carries over to the more general version of the problem with \emph{and+or}-constraints. Conversely, certain $or$-constrained problems can be solved in polynomial time, while they are $\mathcal{NP}$-hard in the $and$-constrained version (Berit \cite{johannes2005complexity}). However, Theorem \ref{thm:complexity_two_machines} shows that it is inaccurate to think that the hardness of \emph{and+or}-constrained problems stems only from the $and$-constraints: While minimizing the makespan or sum of completion times for unit jobs on two machines is in $\mathcal{P}$, both in the $or$- and the $and$-constrained version (Coffman and Graham \cite{coffman1972optimal}, Berit \cite{johannes2005complexity}), we show that scheduling with \emph{and+or}-constraints on just two machines is $\mathcal{NP}$-hard. We suggest studying parameters that measure the ratio between $or$- and $and$-constraints. In Theorem \ref{thm:complexity_two_machines}, we show that $P2 \vert $\emph{and+or-prec}$, p_j=1 \vert C_{max}$ is para-$\mathcal{NP}$-hard parameterized by the number of $and$-successors, while Theorem \ref{thm: konstNach-and+or} proves the same problem to be FPT parameterized by the total number of successors. Open questions in the same spirit include whether restricting only the number of $and$-predecessors or the number of $or$-successors on two machines leads to tractability or retains hardness.\\
Moreover, our results show that it is fruitful to explore parameters that limit job dependencies, especially for problems that become hard due to precedence constraints. This perspective contrasts with the focus on the width in classical $and$-constrained problems, which instead promotes high interrelatedness among jobs. In $and$-constrained problems, parameters of the precedence graph, which reduce interrelatedness, may include the size of a largest clique or tournament and the largest out- or in-degree. In particular, we propose a hybrid approach: combining parameters that limit the portion of input influenced by precedence constraints with parameters that enforce high interrelatedness within affected jobs. 
\bibliographystyle{plain}
\bibliography{lipics-v2021-sample-article}

@article{knop2018scheduling,
  title={Scheduling meets n-fold integer programming},
  author={Knop, Du{\v{s}}an and Kouteck{\`y}, Martin},
  journal={Journal of Scheduling},
  volume={21},
  pages={493--503},
  year={2018},
  publisher={Springer}
}

@article{mnich2015scheduling,
  title={Scheduling and fixed-parameter tractability},
  author={Mnich, Matthias and Wiese, Andreas},
  journal={Mathematical Programming},
  volume={154},
  number={1},
  pages={533--562},
  year={2015},
  publisher={Springer}
}

@InProceedings{nederlof_et_al:LIPIcs.IPEC.2020.25,
  author =	{Nederlof, Jesper and Swennenhuis, C\'{e}line M. F.},
  title =	{{On the Fine-Grained Parameterized Complexity of Partial Scheduling to Minimize the Makespan}},
  booktitle =	{15th International Symposium on Parameterized and Exact Computation (IPEC 2020)},
  pages =	{25:1--25:17},
  series =	{Leibniz International Proceedings in Informatics (LIPIcs)},
  ISBN =	{978-3-95977-172-6},
  ISSN =	{1868-8969},
  year =	{2020},
  volume =	{180},
  editor =	{Cao, Yixin and Pilipczuk, Marcin},
  publisher =	{Schloss Dagstuhl -- Leibniz-Zentrum f{\"u}r Informatik},
  address =	{Dagstuhl, Germany},
  URL =		{https://drops.dagstuhl.de/entities/document/10.4230/LIPIcs.IPEC.2020.25},
  URN =		{urn:nbn:de:0030-drops-133287},
  doi =		{10.4230/LIPIcs.IPEC.2020.25}
}

@inproceedings{goldwasser1996complexity,
  title={Complexity measures for assembly sequences},
  author={Goldwasser, Michael and Latombe, J-C and Motwani, Rajeev},
  booktitle={Proceedings of IEEE International Conference on Robotics and Automation},
  volume={2},
  pages={1851--1857},
  year={1996},
  organization={IEEE}
}

@article{gillies1995scheduling,
  title={Scheduling tasks with AND/OR precedence constraints},
  author={Gillies, Donald W and Liu, Jane W-S},
  journal={SIAM Journal on Computing},
  volume={24},
  number={4},
  pages={797--810},
  year={1995},
  publisher={SIAM}
}

@article{prot2018survey,
  title={A survey on how the structure of precedence constraints may change the complexity class of scheduling problems},
  author={Prot, Damien and Bellenguez-Morineau, Odile},
  journal={Journal of Scheduling},
  volume={21},
  number={1},
  pages={3--16},
  year={2018},
  publisher={Springer}
}

@article{lenstra1978complexity,
  title={Complexity of scheduling under precedence constraints},
  author={Lenstra, Jan Karel and Rinnooy Kan, AHG},
  journal={Operations Research},
  volume={26},
  number={1},
  pages={22--35},
  year={1978},
  publisher={INFORMS}
}

@article{sethi1977complexity,
  title={On the complexity of mean flow time scheduling},
  author={Sethi, Ravi},
  journal={Mathematics of Operations Research},
  volume={2},
  number={4},
  pages={320--330},
  year={1977},
  publisher={INFORMS}
}

@article{mnich2018parameterized,
  title={Parameterized complexity of machine scheduling: 15 open problems},
  author={Mnich, Matthias and Van Bevern, Ren{\'e}},
  journal={Computers \& Operations Research},
  volume={100},
  pages={254--261},
  year={2018},
  publisher={Elsevier}
}

@inproceedings{van2016precedence,
  title={Precedence-constrained scheduling problems parameterized by partial order width},
  author={Van Bevern, Ren{\'e} and Bredereck, Robert and Bulteau, Laurent and Komusiewicz, Christian and Talmon, Nimrod and Woeginger, Gerhard J},
  booktitle={International conference on discrete optimization and operations research},
  pages={105--120},
  year={2016},
  organization={Springer}
}

@article{mohring2004scheduling,
  title={Scheduling with AND/OR precedence constraints},
  author={M{\"o}hring, Rolf H and Skutella, Martin and Stork, Frederik},
  journal={SIAM Journal on Computing},
  volume={33},
  number={2},
  pages={393--415},
  year={2004},
  publisher={SIAM}
}

@article{lee2012flexible,
  title={Flexible job-shop scheduling problems with ‘AND’/‘OR’precedence constraints},
  author={Lee, Sanghyup and Moon, Ilkyeong and Bae, Hyerim and Kim, Jion},
  journal={International Journal of Production Research},
  volume={50},
  number={7},
  pages={1979--2001},
  year={2012},
  publisher={Taylor \& Francis}
}

@article{johannes2005complexity,
  title={On the complexity of scheduling unit-time jobs with OR-precedence constraints},
  author={Johannes, Berit},
  journal={Operations research letters},
  volume={33},
  number={6},
  pages={587--596},
  year={2005},
  publisher={Elsevier}
}

@inproceedings{erlebach2003scheduling,
  title={Scheduling AND/OR-networks on identical parallel machines},
  author={Erlebach, Thomas and K{\"a}{\"a}b, Vanessa and M{\"o}hring, Rolf H},
  booktitle={International workshop on approximation and online algorithms},
  pages={123--136},
  year={2003},
  organization={Springer}
}

@article{bodlaender1995w,
  title={W[2]-hardness of precedence constrained k-processor scheduling},
  author={Bodlaender, Hans L and Fellows, Michael R},
  journal={Operations Research Letters},
  volume={18},
  number={2},
  pages={93--97},
  year={1995},
  publisher={Elsevier}
}

@incollection{graham1979optimization,
  title={Optimization and approximation in deterministic sequencing and scheduling: a survey},
  author={Graham, Ronald Lewis and Lawler, Eugene Leighton and Lenstra, Jan Karel and Kan, AHG Rinnooy},
  booktitle={Annals of discrete mathematics},
  volume={5},
  pages={287--326},
  year={1979},
  publisher={Elsevier}
}

@article{graham1966bounds,
  title={Bounds for certain multiprocessing anomalies},
  author={Graham, Ronald L},
  journal={Bell system technical journal},
  volume={45},
  number={9},
  pages={1563--1581},
  year={1966},
  publisher={Wiley Online Library}
}

@article{coffman1972optimal,
  title={Optimal scheduling for two-processor systems},
  author={Coffman, Edward G and Graham, Ronald L},
  journal={Acta informatica},
  volume={1},
  pages={200--213},
  year={1972},
  publisher={Springer}
}

@article{bessy2019parameterized,
  title={Parameterized complexity of a coupled-task scheduling problem},
  author={Bessy, St{\'e}phane and Giroudeau, Rodolphe},
  journal={Journal of Scheduling},
  volume={22},
  number={3},
  pages={305--313},
  year={2019},
  publisher={Springer}
}

@inproceedings{chen2017parameterized,
  title={Parameterized and approximation results for scheduling with a low rank processing time matrix},
  author={Chen, Lin and Marx, D{\'a}niel and Ye, Deshi and Zhang, Guochuan},
  booktitle={34th Symposium on Theoretical Aspects of Computer Science (STACS 2017)},
  pages={22--1},
  year={2017},
  organization={Schloss Dagstuhl--Leibniz-Zentrum fuer Informatik}
}

@article{happach2021makespan,
  title={Makespan minimization with OR-precedence constraints},
  author={Happach, Felix},
  journal={Journal of Scheduling},
  volume={24},
  number={3},
  pages={319--328},
  year={2021},
  publisher={Springer}
}
\newpage
\appendix

\section{Pseudocodes}
	
	\begin{algorithm}
		\caption{{\sc alg-pre} for scheduling unit time jobs with $k_p$ predecessors}\label{alg:konstVorg}
		\begin{algorithmic}
			\Require An instance of $P \vert$\emph{gen-prec}$, p_j=1 \vert \gamma$ with $\gamma \in \{C_{max}, \sum_j C_j, \sum_j w_j C_j\}$ and predecessors $J_p=\{{i_1},\dots,{i_{k_p}}\} \subseteq J$
			\Ensure Completion times of an optimal schedule\\
			\State Initialize $OPT=\infty$ and $C_j^*=\infty$ for $j \in J$
			\If{$\gamma \in \{C_{max}, \sum_j C_j\}$}
			\State $\mathcal{C}:=\{1,\ldots,k_p\}^{k_p}$
			\Else
            \State $\mathcal{C}:=\{1,\ldots,n\}^{k_p}$
			\EndIf
			\State Let $j_1, \ldots, j_{n-k_p}$ be an ordering of $J\setminus J_p$ such that $w(j_1) \geq \ldots \geq w(j_{n-k_p})$
			\For {$(t_1,\ldots,t_{k_p})\in \mathcal{C}$} 
			\State Define $\mathcal{S}$ by setting $C_{i_\ell}=t_\ell$ for $\ell\in \{1,\ldots,k_p\}$ and set $C_j=\infty$ for $j \in J\setminus J_p$
			
			\If {$(C_{i})_{i \in J_p}$ defines a feasible schedule for $J_p$}
			\For{$r=1, \ldots, {n-k_p}$} 
			
			\State Let $t\in \mathbb{N}_0$ be the earliest time where $j_r$ is available and $\vert \{i \in J\vert C_i=t+1\} \vert <m$
			\State Set $C_j=t+1$
			
			\EndFor
		 
			\If{$\gamma(\mathcal{S})<OPT$ }
			\State Set $C^*_j=C_j$ for $j \in J$ and $OPT=\gamma(\mathcal{S})$ 
			
			\EndIf
			\EndIf
			
			\EndFor

			\State \textbf{return} $(C^*_j)_{j \in J}$

		\end{algorithmic}
	\end{algorithm}
\pagebreak
    \begin{algorithm}
		\caption{{\sc and-alg-succ} for scheduling unit time jobs with $k_s$ successors and $and$-constraints}\label{alg:konstNach}
		\begin{algorithmic}
			\Require An instance of $P\vert prec, p_j=1 \vert \gamma$ with  and successors $J_s=\{{i_1},\dots,{i_{k_s}}\} \subseteq J$
			\Ensure Completion times of an optimal schedule\\
			\State Let $i_1, \ldots, i_{k_s}$ be an arbitrary ordering of the successors
			\State Initialize $OPT=\infty$ and $C_j^*=\infty$ for $j \in J$
			\State Initialize $J_r=J\setminus J_s$ and $\mathcal{C}=\emptyset$
			\For{$(t_1,\ldots,t_{k_s})\in \{1,\ldots,k_s\}^{k_s}$}
			\If {$t_k <t_\ell$ if $i_k$ is a predecessor of $i_\ell$, $k,\ell\in \{1,\ldots,k_s\}$}
			\If{$\vert \{\ell \vert 1 \leq \ell \leq k_s, t_\ell = s\} \vert \leq m$ for all $s\in \{1,\ldots,k_s\}$}
			
			\State $\mathcal{C}=\mathcal{C}\cup \{(t_1+t,\ldots,t_{k_p}+t) \vert t=0, \ldots, n-k_s\}$
			\EndIf
			\EndIf
			\EndFor
			\For{$(t_1,\ldots,t_{k_s})\in \mathcal{C}$}
			\State Define $\mathcal{S}$ by setting $C_{i_\ell}=t_\ell$ for $\ell\in \{1,\ldots,k_s\}$ and set $C_j=\infty$ for $j \in J\setminus J_s$
			
			\For{$q=t_1, \ldots, t_{k_s}$} 
			\For{each predecessor $j$ of a job $k \in J_s$ with $C_k=q$}
			\State Let $t\in \mathbb{N}_0$ be the earliest time where $\vert \{i \in J\vert C_i=t+1\} \vert <m$
			\If{$t+1<C_{i_\ell}$ and $j \in J_r$}
			\State Set $C_j=t+1$
			\EndIf
			\State Set $J_r=J_r \setminus \{j\}$
			\EndFor
			\EndFor
			\For{$j \in J_r$}
			\State Let $t\in \mathbb{N}_0$ be the earliest time where $\vert \{i\in J\vert C_i=t+1\} \vert <m$
			
			\State Set $C_j=t+1$
			\EndFor
			\If{$\gamma(\mathcal{S})<OPT$ }
			\State Set $C^*_j=C_j$ for $j \in J$ and $OPT=\gamma(\mathcal{S})$ 
			
			\EndIf
		
			\EndFor

			\State \textbf{return} $(C^*_j)_{j \in J}$

		\end{algorithmic}
	\end{algorithm}
\pagebreak
\begin{algorithm}
    \caption{{\sc and+or-alg-succ} for scheduling with $k_s$ successors and and+or-constraints}\label{alg:konstNach_andor}
    \begin{algorithmic}
        \Require An instance of $P\vert and+or, p_j=1 \vert \gamma$ with $\gamma \in \{C_{max}, \sum_j C_j\}$, successors $J_s=J_s^{or}\dot\cup J_s^{and}=\{{i_1},\dots,{i_{k_s}}\} \subseteq J$ 
        \Ensure Completion times of an optimal schedule \\
        \State Let $i_1, \ldots, i_{k_s}$ be an arbitrary ordering of the successors
        \State Initialize $OPT=\infty$ and $C_j^*=\infty$ for $j \in J$
        \For{$(t_1,\ldots,t_{k_s})\in \{1,\ldots,k_s\}^{k_s}$}
        
        \If {$t_k <t_\ell$ if $i_k$ is a predecessor of $i_\ell \in J_s^{and}$, $k,\ell\in \{1,\ldots,k_s\}$}
        \If{$\vert \{\ell \vert 1 \leq \ell \leq k_s, t_\ell = s\} \vert \leq m$ for all $s\in \{1,\ldots,k_s\}$}
        
        \State $\mathcal{C}=\mathcal{C}\cup \{(t_1+t,\ldots,t_{k_p}+t) \vert t=0, \ldots, k_s\}$
        \EndIf
        \EndIf
        \EndFor
        \For{$(t_1,\ldots,t_{k_s})\in \mathcal{C}$}
        \State Define a schedule $\mathcal{S}$ by setting $C_{i_\ell}=t_\ell$ for $\ell\in \{1,\ldots,k_s\}$ and set $C_j=\infty$ for $j \in J\setminus J_s$
            \State Initialize $I_s=\emptyset$
        \For{$j \in J_s^{or}$}
        \If{some $p \in J_s$ is a predecessor of $j$ with $C_p<C_j$}
        \State Set $\psi_j=x_p$
        \State Add $j$ to $I_s$
        \EndIf
        \EndFor
        \State Sort $j_1, \ldots ,j_r \in J_s^{or}$ in order of non-decreasing completion time
        \For{$(S_j)_{j \in J_s^{or}} \in \prod_{j \in J_s^{or}} \{S \subseteq J_s \vert j \in S\}$}
            \State Initialize $I^{or}=J_s^{or}\setminus I_s$
            \For{$\ell=1, \ldots ,r$ }
            \If{$j_\ell \in I^{or}$}
            \If{there exists $i \in J$ such that $i$ is a predecessor of all in jobs in $S_{j_\ell}$ and is not a predecessor of any job in $(\{j_\ell, \ldots, j_r\}\cup J_s^{and})\setminus S_{j_\ell}$}
            
            \For{$s=\ell, \ldots, r$}
            \If{$i$ is a predecessor of $j_s$}
            \State Set $\psi_{j_s}=x_i$
            \State Remove $j_s$ from $I^{or}$
            \EndIf
            \EndFor
            \Else{Go to the next selection of $(S_j)_{j \in J_s^{or}}$}

            \EndIf
            \EndIf
            \EndFor

            \State Let $\mathcal{S}=(C_j)_{j\in J}$ be the completion times that {\sc and-alg-succ} computes for $(J,(\psi_j)_j)$ and the choice of $(C_j)_{j \in J_s}$

			\If{$\gamma(\mathcal{S})<OPT$ }
			\State Set $C^*_j=C_j$ for $j \in J$ and $OPT=\gamma(\mathcal{S})$ 
			
			\EndIf
		
			\EndFor
            \EndFor

			\State \textbf{return} $(C^*_j)_{j \in J}$

    \end{algorithmic}
\end{algorithm}
\pagebreak
    \section{Proofs}

    \subsection{Constant Predecessors}\label{app:konstVorg}
	\konstVorg*
		
	\begin{proof}
		Clearly, {\sc alg-pre} computes a feasible schedule. Consider an optimal schedule $\mathcal{S}^*$ with completion times $C_j^*$, $j\in J$ such that $(C^*_{i_1}, \ldots, C^*_{i_{k_p}}) \in \mathcal{C}$. (For $\gamma\in  \{C_{max}, \sum_j C_j\}$, the existence of such $\mathcal{S}^*$ follows from Lemma \ref{lem: predecessors_first}). Consider the iteration of {\sc alg-pre} where $C_j=C_j^*$ for all $j \in J_p$. Clearly, the algorithm computes the completion times $(C_j)_{j \in J}$ of a feasible schedule $\mathcal{S}$. Denote by $k_p^t$ the number of predecessor jobs that $\mathcal{S}$ and $\mathcal{S}^*$ process during time slot $t\in \{1, \ldots, n\}$.\\
        Assume the jobs in $J \setminus J_p$ have $1\leq q\leq n-k$ different weights, that we denote by $w^1 \geq w^2 \geq \ldots \geq w^q$. We consider the vector $z(t)=(z_1(t), \ldots, z_q(t)) \in \mathbb{N}^q$ where $z_i(t)$, $i \in \{1, \ldots, q\}$, denotes the number of jobs in $J \setminus J_p$ with weight $w^i$ that $\mathcal{S}$ processes during time slot $t\in  \mathbb{N}$. For $\mathcal{S}^*$, we define $z^*(t)$, $t\in  \mathbb{N}$, analogously. If $z(t)=z^*(t)$ for all $t \in \mathbb{N}$, then $(C_{max}, \sum_j C_j, \sum_j w_jC_j)=(C^*_{max}, \sum_j C^*_j, \sum_j w_j C^*_j)$ and $\mathcal{S}$ is optimal. Otherwise, let $T:=\min\{t \vert z(t) \neq z^*(t)\}$ be the smallest time slot for which $z(T),z^*(T)$ are not equal. \\
		Assume first that $z(T)$ is lexicographically smaller than $z^*(T)$. Let $i\in \{1,\ldots,q\}$ be such that $z_s^*(T)=z_s(T)$ for $s<i$ and $z_i^*(T)>z_i(T)$. It holds that $\sum_{t=1}^T z_i(t)<\sum_{t=1}^Tz_i^*(t)$. Subsequently, there exists a job $j \in J$ with weight $w_j=w^i$, such that $C^*_j \leq T$ but $C_j > T$. Observe that $j$ is no predecessor, otherwise $C_j=C^*_j$. Since $C^*_j \leq T$, job $j$ is available at time $T-1$ in $\mathcal{S}^*$ and thus also in $\mathcal{S}$. Since {\sc alg-pre} sets $C_j > T$, the following holds: By the iteration $r\in\{ 1, \ldots, {n-k_p}\}$ in which $j=j_r$ is considered, the algorithm must have already occupied all machines during time slot $T$, namely with predecessor jobs or jobs with weight at least $w_j=w^i$. Then
		\begin{align*}
			m=k_p^T+\sum_sz_s(T)&=k_p^T+\sum_{s<i}z_s(T)+z_i(T)=k_p^T+\sum_{s<i}z^*_s(T)+z_i(T) \\&<k_p^T+\sum_{s<i}z^*_s(T)+z^*_i(T) \leq m, 
		\end{align*}
		contradiction. \\		
		Assume now that $z^*(T)$ is lexicographically smaller than $z(T)$ and let $i\in \{1,\ldots,q\}$ be such that $z_s^*(T)=z_s(T)$ for $s<i$ and $z_i^*(T)<z_i(T)$. By analogous arguments, there exists a job $j \in J\setminus J_p$ with weight $w^i$ such that $C_j \leq T$, $C^*_j > T$ and $j$ is available at time $T-1$ in $\mathcal{S}^*$. Assume first $\gamma=C_{max}$, in which case we have uniform weights and $i=q=1$. Since $z_i^*(T)<z_i(T)\leq m$, there exists an idle machine in $\mathcal{S}^*$ during time slot $T$. We alter the optimal schedule $S^*$ by setting $C^*_j=T$. Clearly, this does not negatively affect the feasibility nor the makespan of $\mathcal{S}^*$. Repeat until $z^*(t)=z(t)$ for all $t \in \mathbb{N}$. Then $\mathcal{S}^*$ and $\mathcal{S}$ have the same makespan and thus $\mathcal{S}$ is optimal. Assume now $\gamma \in \{\sum_j C_j, \sum_j w_jC_j\}$. The optimal schedule $\mathcal{S}^*$ occupies all machines during time slot $T$ with predecessor jobs or jobs of weight at least $w_j=w^i$: Otherwise, swapping $j$ with a job of smaller weight or placing it on an idle machine during time slot $T$ would yield a feasible solution with smaller (weighted) sum of completion time. As before,
		
			$$m=k_p^T+\sum_{s<i}z^*_s(T)+z^*_i(T)<k_p^T+\sum_{s<i}z_s(T)+z_i(T)\leq m$$
		
		\noindent yields a contradiction.\\
		{We conclude with a brief runtime analysis: We sort all non-predecessor jobs with respect to weight once ($\mathcal{O}(n\cdot log(n))$ time). Then the algorithm checks $\vert \mathcal{C} \vert$ configurations $(t_1,\ldots,t_{k_p})\in \mathcal{C}$ of predecessor completion times. For each such configuration, the algorithm inserts $\mathcal{O}(n)$ jobs $j\in J\setminus J_p$. For each job $j$, it checks for $\mathcal{O}(n)$ time steps $t$ whether $j$ is available at time $t$. For any such time step $t$, the algorithm searches for an idle machine among $m$ machines. Since, $\vert \mathcal{C} \vert \in \mathcal{O}(k^k)$ if $\gamma\in  \{C_{max}, \sum_j C_j\}$ and $\vert \mathcal{C} \vert \in \mathcal{O}(n^k)$ otherwise, the claimed asymptotic runtime follows.}		
	\end{proof}
\subsection{Parameterization by the number of successors for $and$-constraints}\label{app:konstNachfolgand}
	
	\konstNachand*
\begin{proof}
   By Lemma \ref{lem: successors_span}, there exists an optimal schedule $\mathcal{S}^*$ with completion times $C_j^*$, $j\in J$ such that $(C^*_{i_1}, \ldots, C^*_{i_{k_p}}) \in \mathcal{C}$. Consider the iteration of {\sc and-alg-succ} where $C_j=C_j^*$ for all $j \in J_s$. 
   Since $\mathcal{S}^*$ is a feasible schedule with identical successor completion times, there is enough space on the machines prior to each successor $j \in J_s$ for {\sc and-alg-succ} to schedule all of $j$'s predecessors. Hence, $(C^*_{i_1}, \ldots, C^*_{i_{k_p}}) $ is not discarded and {\sc and-alg-succ} computes a feasible schedule $\mathcal{S}$. 
   We proceed by proving that $\mathcal{S}$ is in fact optimal. We define $z(t)$ as the number of jobs in $\{j \in J \setminus J_s\vert C_j=t\}$, which $\mathcal{S}$ processes during time slot $t\in  \mathbb{N}$. For $\mathcal{S}^*$, we define $z^*(t)$, $t\in  \mathbb{N}$, analogously. If $z(t)=z^*(t)$ for all $t \in \mathbb{N}$, then $(C_{max}, \sum_j C_j)=(C^*_{max}, \sum_j C^*_j)$ and $\mathcal{S}$ is optimal. Otherwise, let $T:=\min\{t \vert z(t) \neq z^*(t)\}$ be the smallest time slot for which $z(T),z^*(T)$ are not equal. \\
Assume first that $z(T)<z^*(T)$. It holds that $\sum_{t=1}^T z(t)<\sum_{t=1}^Tz^*(t)$. Subsequently, there exists a job $j \in J$, such that $C^*_j \leq T$ but $C_j > T$. Observe that $j$ is no successor, otherwise $C_j=C^*_j$. When {\sc and-alg-succ} assigns a completion time to $j$, there is an idle machine during time slot $T$ because $z(T)<z^*(T)\leq m$ and {\sc and-alg-succ} assigns a completion time $C_j \leq T$ to $j$, contradiction. If $z(T)>z^*(T)$, there exists a job $j \in J\setminus J_s$, such that $C^*_j > T$ but $C_j \leq T$. The optimal schedule $\mathcal{S}^*$ has an idle machine during time slot $T$ because $z^*(T)<z(T)\leq m$. Adapting the optimal solution such that $C^*_j=T$ does not negatively affect neither the optimality of $\mathcal{S}^*$ nor its feasibility since $j$ has no predecessors. We repeat this argument until $z(t)=z^*(t)$ for all $t \in \mathbb{N}$, which proves that $\mathcal{S}$ is optimal. We omit a runtime analysis since it is analogous to the proof of Theorem \ref{alg:konstVorg}.
\end{proof}

\subsection{Parameterization by the number of successors for $and+or$-constraints}\label{app:konstNachfolgand+or}

\represent*
\begin{proof}
	Assume that $L$ is infeasible and let $C=(j_1, \ldots, j_r)$ be a cycle in the precedence graph. We define $j_0:=j_r$ and consider a job $j_\ell$ in the cycle, $\ell \in \{1, \ldots, r\}$. If $j_\ell \in J_s^{and}$ is an $and$-job in $I$, then $j_{\ell-1}$ is also a predecessor of $j_\ell$ in $L'$. If $j_\ell \in J_s^{or}$, then $j_{\ell-1}=i_j$. Since $j_{\ell-1}$ is a successor in $L$, it is also a successor in $I$ and thus $j_{\ell-1}=i_j=i'_j$ is also a predecessor of $j_\ell$ in $L'$. Therefore, the precedence graph of $L'$ contains the same cycle $C$, due to which $L'$ is infeasible.\\
	Assume now that $L$ is feasible, $\sigma=\sigma'$ and $i_j=i'_k$ whenever $i'_j=i_k$, $j,k \in J_s^{or}$. Let $(C_j)_{j \in J}$ be the completion times of an optimal schedule $\mathcal{S}$. We define a schedule $\mathcal{S'}$ with completion times $(C'_j)_{j \in J}$ for $L'$. We denote the sets in $\{\sigma(j)\cap J_s^{or}\vert j \in J_s^{or}\}$ by $S_1, \ldots, S_r$. Observe that $(S_\ell)_{\ell=1}^r$ is a partition of $J_s^{or}$ and that $S_\ell =\sigma(j)$ precisely if $S_\ell=\sigma'(j)$, $j\in J_s^{or}, \ell =1, \ldots, r$. For each $\ell=1, \ldots, r$, let $i(\ell),i'(\ell) \in J_s^{or}$ be the common predecessor with $i(\ell)=i_k, i'(\ell)=i'_k$ for all $k\in S_\ell$. We obtain $\mathcal{S}'$ from $\mathcal{S}$ by swapping the positions of $i(\ell)$ and $i'(\ell)$ for all $\ell =1, \ldots,r$. 
	Observe that this is in fact the same as setting $C'_{i'_k}=C_{i_k}$ and $C'_{i_k}=C_{i'_k}$ for $k \in J_s^{or}$. All other jobs $j\in J$ remain as in $\mathcal{S}$, i.e. $C'_j=C_j$. Note that the schedule $\mathcal{S'}$ is well-defined, since we require $i'(\ell)=i(p)$ if $i(\ell)=i'(p)$ for $\ell, p =1, \ldots, r$. 
	Clearly, $\gamma(\mathcal{S}) =\gamma(\mathcal{S}') $. It remains to show that $\mathcal{S}'$ is feasible. Let $k \in J_s$ be a successor job in $I$. Observe that $C_k=C'_k$: If $k \neq i_\ell, i'_\ell$ for all $\ell \in J_s^{or}$, then $C'_k=C_k$ by definition of $C_k$. Otherwise, if $k=i_\ell$ or $k=i'_\ell$ for some $\ell \in J_s^{or}$, then $\ell \in I_s$ and thus $k=i_\ell=i'_\ell$ and $C'_{k}=C_{k}$. Assume now that $k \in J_s^{and}$ is an $and$-job in $I$ and let $j \in J$ be a predecessor of $k$ in $L'$. Clearly, $j$ is also a predecessor of $k$ in $L$ and thus $C_j < C_k$. If $C'_{j}=C_{j}$, then $j$ precedes $k$ in $\mathcal{S}'$ because $C'_{j}=C_{j}<C_{k}=C'_{k}$. Otherwise, $j=i_\ell$ or $j=i'_\ell$ for some $\ell \in J_s^{or}$ with $i_\ell\neq i'_\ell$. Assume first that $j=i_\ell$. Since $k$ is an $and$-job in $I$ with predecessor $i_\ell=j$, it follows that $k \in \sigma(\ell)=\sigma'(\ell)$. By definition of $\sigma'$, job $i'_\ell$ is a predecessor of $k$ in $I$ and $L$ and thus $C'_{j}=C_{i'_\ell}<C_k=C'_k$. We can show analogously that $C'_{j}<C'_k$ if $j=i'_\ell$. 
	Finally, let $k \in J_s^{or}$ be a former $or$-job. By definition of $L'$, the job $i'_k$ is the only predecessor of $k$ in $L'$ and $i_k$ is the only predecessor of $k$ in $L$. Subsequently, $C'_{i'_k}=C_{i_k}<C_k=C'_k$, which completes the proof.
\end{proof}

	\konstNachandplor*

\begin{proof}
	Let $J_s=J_s^{or} \dot \cup J_s^{and}$ be the set of successors with $\vert J_s \vert =k_s$ and $q:=\vert J_s^{or} \vert$. By Lemma \ref{lem: successors_span}, there exists an optimal schedule $\mathcal{S}^*$ with completion times $C_j^*$, $j\in J$ such that $(C^*_{i_1}, \ldots, C^*_{i_{k_p}}) \in \mathcal{C}$. Consider the iteration of {\sc and+or-alg-succ} where $C_j=C_j^*$ for all $j \in J_s$ and let $\mathcal{S}$ denote the schedule defined by $(C_j)_j$. Let $j_1, \ldots, j_q$ be the ordering, in which {\sc and+or-alg-succ} iterates over the $or$-jobs. Moreover, we denote by $I_s$ the set of $or$-jobs $j \in J_s^{or}$, to which {\sc and+or-alg-succ} assigns a predecessor $i_j \in J_s$. \\
	The idea of the proof is the following: First, we pick $(i^*_j)_{j \in J_s^{or}}$ such that $i^*_j$ is a predecessor of $j$ with $C^*_{i^*_j} < C^*_j$, $j \in J_s^{or}$. We argue, that for a certain choice of $(S_j)_j$, the algorithm computes representatives $(i_j)_{j \in J_s^{or}}$ such that $(i_j)_j,(i^*_j)_j$ as well as the associated common predecessor functions $\sigma, \sigma^*$ comply with the conditions in Lemma \ref{lem: representatives}. \\
	More precisely, we pick $(i^*_j)_{j \in J_s^{or}}$ as follows: Consider the $or$-jobs $j$ in order $j=j_1, \ldots, j_q$. We set $i^*_j=i_j$ for $j \in I_s$. If $j=j_i$ and no $i^*_j$ has been chosen yet, we pick an arbitrary predecessor $i^*_j$ with $C^*_{i^*_j} < C^*_j$ but set $i^*_k=i^*_j$ for all $k=j_{i+1}, \ldots, j_q$, of which $i^*_j$ is also a predecessor. We denote by $\sigma^*:J_s^{or} \rightarrow 2^{J_s}$ the common predecessor function of $I$ w.r.t. $(i^*_j)_{j \in J_s^{or}}$ and consider the iteration where {\sc and+or-alg-succ} sets $S_{j_i}=\sigma^*(j_i)\setminus\{j_1, \ldots, j_{i-1}\}$ for $i=1, \ldots, q$. Observe that whenever {\sc and+or-alg-succ} is required to pick $i_j$ according to $S_j$, then $S_j=\sigma^*(j)$, otherwise $i_j$ would have been chosen in an earlier iteration. Now we argue that {\sc and+or-alg-succ} can find an $i_j$ for $j=j_r$ such that 1. $i_j$ is a predecessor of all jobs in $S_j$ and 2. $i_j$ is not a predecessor of any job in $(J_s^{and}\cup \{j_{r+1}, \ldots, j_q\})\setminus S_j$. Clearly, $i^*_j$ is a common predecessor of all jobs in $S_j=\sigma^*(j)$. If $k\in J_s^{and}\setminus \sigma^*(j)$ and $i^*_j$ is a predecessor of $k$, then $k\in \sigma^*(j)$ by the definition of $\sigma^*$, contradiction. If $k\in \{j_{r+1}, \ldots, j_q\}\setminus \sigma^*(j)$ and $i^*_j$ is a predecessor of $k$, then $i^*_k=i^*_j$ by the choice of $i^*_k$, and thus $k \in \sigma^*(j)$, contradiction. We denote by $\sigma:J_s^{or} \rightarrow 2^{J_s}$ the common predecessor function of $I$ w.r.t. $(i_j)_{j \in J_s^{or}}$.\\ 
	Now we prove that $(i_j)_j,(i^*_j)_j$ and $\sigma, \sigma^*$ comply with the conditions in Lemma \ref{lem: representatives}. By the choice of $(i^*_j)_j$, it is easy to see that $I_s:=\{j \in J_s^{or} \vert i_j \in J_s\}=\{j \in J_s^{or} \vert i^*_j \in J_s\}$, $i_j=i^*_j$ for all $j \in I_s$. Before we check the other conditions, we show that $i_j$ is a predecessor of all jobs in $\sigma^*(j)$ while $i_j$ is not a predecessor of any job in $(J_s^{and}\cup \{j_{r+1}, \ldots, j_q\})\setminus \sigma^*(j)$, $j=j_r \in J_s^{or}$. The claim obviously holds if the algorithm picks $i_j$ according to $S_j$. If $j \in I_s$, the claim follows from $i_j=i^*_j$ and the choice of $(i^*_\ell)_\ell$. Otherwise, $i_j=i_{j_p}$ for some $p<r$, where $i_{j_p}$ was picked according to $S_{j_p}$. Then all successors of $i_{j_p}$ in $\{j_{p+1}, \ldots, j_q\}$ are in $S_{j_p}=\sigma^*(j_p)$. Therefore, $j \in \sigma^*(j_p)=\sigma^*(j)$ and $i_j=i_{j_p}$ is a predecessor of all jobs in $S_{j_p}=\sigma^*(j_p)=\sigma^*(j).$ Moreover, all successor of $i_j=i_{j_p}$ lies in $(J_s^{and}\cup \{j_{r+1}, \ldots, j_q\})\setminus \sigma^*(j)\subset (J_s^{and}\cup \{j_{p+1}, \ldots, j_q\})\setminus \sigma^*(j_p)$. We proceed by showing that $\sigma=\sigma^*$. Let first $j \in I_s$ and consider $k \in \sigma^*(j) \cup \sigma(j)$. If $k \in J_s^{or}$ then $i_k=i^*_k=i^*_j=i_j$, where the first equality is implied by $k \in I_s$. Hence, $k \in \sigma^*(j) \cap \sigma(j)$. If $k \in J_s^{and}$ then $i^*_j=i_j$ is a predecessor of $k$ which implies again $k \in  \sigma^*(j) \cap  \sigma(j)$. Assume now $j \in J_s^{or}\setminus I_s$ and let $k \in \sigma^*(j)\cap J_s^{and}$. Since $i_j$ is a predecessor of all jobs in $S_j=\sigma^*(j)$, $i_j$ is a predecessor of $k$ and thus $k\in \sigma(j)$. Let now $k \in \sigma^*(j)\cap J_s^{or}$. Then $\sigma^*(k)=\sigma^*(j)$ and the algorithm sets $i_j=i_k$. Therefore, $k\in \sigma(j)$. We now show that $\sigma(j)\subseteq \sigma^*(j)$ for all $j \in J_s^{or}\setminus I_s$. Let $j \in J_s^{or}\setminus I_s$ and consider first $k \in \sigma(j)\cap J_s^{and}$. Then $k$ is a successor of $i_j$ and, as seen above, all successors of $i_j$ in $J_s^{and}$ lie in $\sigma^*(k)$. Consider now $k \in \sigma(j)\cap J_s^{or}$, which implies that $i_k=i_j$. Let w.l.o.g. $j=j_p$ and $k=j_r$ with $p<r$. Then all successors of $i_{j}$ in ${j_{p+1}, \ldots, j_q}$, including $k$, lie in $\sigma^*(j)$. To complete the proof, we need to show that $i^*_k=i_\ell$ implies $i_k=i^*_\ell$ for all $k,\ell \in J_s^{or}$. Consider $k,\ell \in J_s^{or}$ with $i^*_k=i_\ell$ and assume first that $k=j_p$, $\ell=j_r$ with $p<r$. Then all successors of $i^*_k$ in $\{j_{p}, \ldots, j_q\}$, including $\ell$, lie in $\sigma^*(k)=\sigma(k)$ and thus $i^*_\ell=i^*_k=i_\ell=i_k$. Analogously, we can show that $i^*_\ell=i^*_k=i_\ell=i_k$ if $\ell=j_p$, $k=j_r$ with $p<r$. \\
	Finally, we briefly discuss the runtime. For each $or$-job $j$, we have $\mathcal{O}(2^{k_s})$ possibilities to choose a set $S_j$, which results in $\mathcal{O}(2^{k_s^2})$ many possibilities for choosing $(S_j)_{j \in J_s^{or}}$. For each configuration, we need to apply Algorithm {\sc and-alg-succ}. The overall runtime follows from Proposition \ref{thm: konstNach-and}.
	\end{proof}

\subsection{Parameterization by the number of successors for $or/and$-constraints}\label{app:konstNachfolgerorand}
\weins*

\begin{proof}
	We prove the claim by an FPT-reduction from Independent Set. Let $I'=(G,k)$ be an instance of Independent Set where $G=(V,E)$ is an undirected graph with $n':=\vert V \vert$ nodes and $k \in \mathbb{N}_0$ is the required size of an independent set, w.l.o.g. $k<n$. By $\overline{\mathcal{N}}(v):=\{w \in V \vert \{v,w\} \in E\} \cup \{v\}$ we denote the closed neighborhood of a vertex $v \in V$ and we define $n'_v:= \vert \overline{\mathcal{N}}(v) \vert$ as the number of vertices in the closed neighborhood of $v$. We will first reduce Independent Set to $P\vert or/and$-$prec, p_j=1 \vert C_{max}$. The proof idea is the following: We construct from $I'$ an instance $I$ of $P\vert or/and$-$prec, p_j=1 \vert C_{max}$ with $n'$ machines and job set $J$ and require a makespan of at most $2(k+1)$. We introduce $k$ jobs $u_1, \ldots, u_k$ where $u_i$ represents the $i$-th node in an independent set, $i=1, \ldots, k$. The precedence constraint of a job $u_i$, $i=1, \ldots,k$, is composed of $n'$ clauses $K_{iv}$, $v\in V$. Each clause $K_{iv}$ is associated with a node $v \in V$ and we interpret $v$ as the $i$-th node in the independent set if $K_{iv}$ is satisfied by the jobs that precede $u_i$, $i=1, \ldots, k$, $v\in V$. We use a chain of $k+1$ auxiliary jobs $h_1,  \ldots, h_{k+1}$ to make sure that all jobs $u_1, \ldots, u_k$ have to be completed by time $k+1$. This leaves a total of at most $k \cdot n'$ time slots to schedule the predecessors of $u_1, \ldots, u_k$ which are needed to make $u_1, \ldots, u_k$ available. However, each clause $K_{iv}$ is composed of $n'+i-1$ literals: 1. $n'_v$ auxiliary jobs $b_{iu}$, $u \in \overline{\mathcal{N}}(v)$, 2. $n'-n'_v$ jobs $N_{iu}$ which represent all nodes $u$ outside the closed neighborhood of $v$ and 3. $i-1$ jobs $N_{jv}$, $j=1, \ldots, i-1$. Observe that the first $n'$ predecessors in $K_{iv}$ cannot be shared predecessors with any other $u_j$, $j\neq i$, and thus each of the last $i-1$ predecessors $N_{jv}$, $j=1, \ldots,i-1$, has to be a common predecessor with $u_j$, due the space restrictions. Subsequently, $v$ is disjoint and non-adjacent to the previous $i-1$ nodes in the independent set. An example of the construction and a feasible schedule is depicted in Figure \ref{fig:w_1}.\\
	Now we describe the construction in more detail. We introduce $k$ jobs $u_i$, $i=1, \ldots, k$, where each $u_i$ represents a node in the independent set. Moreover, we introduce jobs $N_{iv}$, $i=1, \ldots, k$, $v \in V$, where job $N_{iv}$ is interpreted as $v$ lying outside the closed neighborhood of the $i$-th vertex in the independent set. We need $k+1$ auxiliary jobs $h_{1}, \ldots, h_{k+1}$ which will be used to make sure that all jobs $u_1, \ldots, u_k$ have to be completed by time $k+1$ (if we were to meet the required makespan). Finally, we introduce auxiliary jobs $b_{iv}$, $i=1, \ldots, k$, $v\in V$, which will be used in clauses, to make sure that each $u_i$ has precisely $n'$ non-shareable predecessors. The precedence constraints are defined as follows: The auxiliary jobs $h_{1}, \ldots, h_{k+1}$ have to be scheduled in a chain; that is, $\psi_{h_{i+1}}=x_{h_i}$ for all $i=1, \ldots,k$. The first auxiliary job $h_1$ cannot be scheduled unless all $u_i$ have been completed: $\psi_{h_1}=\bigwedge_{i=1}^k u_i$. The precedence constraint of a job $u_i$ is the disjunction $\psi_{u_i}=\bigvee_{v\in V}K_{iv}$ of $n'$ clauses $K_{iv}$, where 
    \begin{equation*}
        K_{iv}:=\bigwedge_{u \notin \overline{\mathcal{N}}(v)} x_{N_{iu}} \wedge \bigwedge_{u \in \overline{\mathcal{N}}(v)} x_{b_{iu}} \wedge \bigwedge_{j\leq i}x_{N_{jv}},\text{ }i=1, \ldots, k, v\in V.
    \end{equation*}
	We prove that $G$ has an independent set of size $k$ precisely if $OPT(I) \leq 2k+2$. Let first $\mathcal{S}$ be a schedule with completion times $(C_j)_{j \in J}$ and $C_{max} \leq 2k+2$. For each $u_i$, $i =1, \ldots, k$, let $v_i\in V$ be a node such that $K_{iv_i}$ is satisfied by the jobs $P_i:=\{j \in J\vert x_j \text{ in }K_{iv_i}\}$ that precede $u_i$ in $\mathcal{S}$. Set $X:=\{v_i\vert i=1, \ldots,k\}$ and show that $X$ is an independent set of size $k$. Since $C_{h_{k+1}}\leq C_{max}\leq 2k+2$ and $C_{h_{r}}< C_{h_{r+1}}$, $r=1, \ldots, k$, we know that $C_{u_i} < C_{h_1} \leq k+2$ for all $i=1, \ldots, k$. Subsequently, we have at most $k \cdot n'$ time slots to schedule all jobs in $\cup_{i=1}^kP_i$. Observe that $S_i:=\{N_{iu}\vert u \notin \overline{\mathcal{N}}(v_i)\}\cup \{x_{b_{iu}} \vert u \in \overline{\mathcal{N}}(v_i)\} \subseteq P_i$ is a subset of predecessors with $\vert S_i \vert = n'$ and $S_i \cap S_j=\emptyset$ for $i,j \in \{1, \ldots, k\}$, $i\neq j$. From $\vert\cup_{i=1}^kP_i\vert\leq k\cdot n'$, it follows that $S_i=P_i$ for $i\in \{1, \ldots, k\}$ and $\{j \in J\vert C_j \leq k\}=\cup_{i=1}^kS_i$.
	Assume now by contradiction that there exist $i,j\in \{1,\ldots,k\}$ with $i< j$ such that $v_j \in \overline{\mathcal{N}}(v_i)$. Then $N_{iv_j}\notin \cup_{r=1}^kS_r$ and hence $C_{N_{iv_j}} >k$, a contradiction to $N_{iv_j} \in P_j$.\\
	Conversely, let $X=\{v_1, \ldots, v_k\}$ be an independent set and prove that $OPT(I)\leq 2k+2$. We set $C_{u_i}=k+1$ for $i=1, \ldots, k$ and we set $C_{h_j}=k+j+1$, $j=1, \ldots,k+1$. For $i\in \{1, \ldots , k\}$,  $u\in V\setminus \overline{\mathcal{N}}(v_i)$ and $w\in  \overline{\mathcal{N}}(v_i)$, we set $C_{N_{iu}}=i=C_{b_{iw}}$. Each of the remaining $n'\cdot k$ jobs is scheduled as early as possible. It is easy to see that the resulting schedule is feasible with makespan $C_{max}=2k+2$. Observe that the constructed schedule processes $2k+2$ jobs on a single machine without idle time and is perfect on the other $n'-1$ machines, due to which the proof extends to $\gamma=\sum_jC_j$. 
\end{proof}

\subsection{Parameterization by the number of machines}\label{app:twomachines}

	\twomachines*
\begin{proof}
	We prove the claim by a reduction from vertex cover. Let $I=(G,k)$ be an instance of vertex cover, where $G=(V,E)$ is an undirected graph and we have to decide whether $G$ contains a vertex cover of size $k$. We assume w.l.o.g. that $k<\vert V \vert$, otherwise $I$ is trivial. We construct from $I$ an instance $I'$ of $P2 \vert $\emph{and+or-prec}$, p_j=1 \vert \gamma$, $\gamma\in \{C_{max}, \sum_j C_j\}$, where the idea is the following: We introduce a job for each edge and for each vertex in $G$. Moreover, we introduce $\ell:=\vert V \vert +  \vert E \vert$ auxiliary jobs to block the second machine. Using precedence constraints, we will force a feasible schedule to schedule all jobs in a vertex cover before finishing the edge jobs. 
    The auxiliary jobs are used to ensure that any perfect schedule can start only up to $k$ vertex jobs prior to finishing the edge jobs. \\More precisely, we set $J:= V \cup E \cup B$ where $B=\{b_1, \ldots, b_{\ell}\}$ is the set of auxiliary jobs. All vertex jobs, as well as the first auxiliary job, may start at any time, i.e. $\psi_j=\mathbbm{1}$ if $j \in V \cup \{b_1\}$. Edges may not be scheduled unless one of their end-vertices has been scheduled already, i.e. $\psi_{e}=x_{u} \lor x_v$ for $e =\{u,v\} \in E$. Auxiliary jobs shall not be scheduled simultaneously, thus, we make $b_{i-1}$ a predecessor of $b_i$ for each $i=2, \ldots, \ell$. Moreover, we make sure that all edge jobs have been scheduled by the time we start the $(k+\vert E \vert+1)$-st auxiliary job to ensure that a perfect schedule completes all edge jobs by the time $t=k+\vert E \vert$. That is, we set $\psi_{b_{k+|E|+1}}=x_{b_{k+|E|}} \land \bigwedge_{e\in E} x_{e}$. For the remaining auxiliary jobs, we set $\psi_{b_{i}}=x_{b_{i-1}} $, $i \in \{2, \ldots, \ell\} \setminus \{k+\vert E\vert+1\}$. The construction is visualized in Figure \ref{fig:constMachines}.
	Observe that a schedule is perfect precisely if it yields $C_{max}=\ell$ or, equivalently, if $\sum_j C_j= \ell\cdot(\ell+1)$. Any perfect schedule processes all $b_i$ uninterruptedly in order of increasing $i=1, \ldots ,\ell$. By time $k+m$, all edge jobs need to be processed because  ${b_{k+1+m}}$ starts at time $k+m$ and each edge is a required predecessor of ${b_{k+1+m}}$. This leaves room for at most $k$ vertex jobs to be scheduled until time $k+m$. The vertices that are scheduled before time $k+m$ yield a vertex cover, otherwise some edge's precedence constraint is violated.
	Let now $X \subseteq V$ be a vertex cover with $\vert X \vert \leq k$. Then it is easy to see that the following schedule is feasible and perfect: Process all auxiliary jobs $b_i$ on one machine in order of increasing $i=1, \ldots, \ell$. The other machine starts by processing all jobs in $X$, continues to schedule all edge jobs $e_i$ in order of increasing $i=1, \ldots, \vert E \vert$, and terminates by scheduling all remaining vertex jobs in arbitrary order. Overall, $I$ is a yes-instance, precisely if there exists a schedule for $I'$ with $C_{max}=\ell$, ($\sum_j C_j= \ell\cdot(\ell+1)$). 
\end{proof}

\end{document}